\documentclass[12pt]{article}
\usepackage[margin=1in]{geometry}
\usepackage[T1]{fontenc}
\usepackage[utf8]{inputenc}
\usepackage{amsmath,amssymb,amsfonts,amsthm,mathrsfs}
\usepackage{color}
\usepackage[hidelinks]{hyperref}

\numberwithin{equation}{section}
\newtheorem{theorem}{Theorem}[section]
\newtheorem{lem}[theorem]{Lemma}
\newtheorem{prop}[theorem]{Proposition}

\newtheorem{ass}[theorem]{Assumption}
\theoremstyle{definition}
\newtheorem{definition}[theorem]{Definition}
\theoremstyle{remark}
\newtheorem{rem}[theorem]{Remark}

\newcommand{\beq}[1]{\begin{equation} \label{#1}}
\newcommand{\eeq}{\end{equation}}
\newcommand{\bea}{\bed\begin{array}{rl}}
\newcommand{\eea}{\end{array}\eed}
\newcommand{\bed}{\begin{displaymath}}
\newcommand{\eed}{\end{displaymath}}
\newcommand{\barray}{\begin{array}{ll}}
\newcommand{\earray}{\end{array}}
\newcommand{\disp}{\displaystyle}
\newcommand{\ad}{&\!\disp}
\newcommand{\aad}{&\disp}
\newcommand{\beqa}[1]{\begin{equation}\label{#1}\barray}
\newcommand{\eeqa}{\earray\end{equation}}
\newcommand{\e}{\varepsilon}

\newcommand{\La}{\Lambda}

\newcommand{\ga}{\gamma}
\newcommand{\Ga}{\Gamma}
\newcommand{\dl}{\delta}
\newcommand{\Dl}{\Delta}
\newcommand{\cd}{(\cdot)}

\def\phi{\varphi}
\def\indi{{\bf 1}}
\def\half{\frac{1}{2}}

\newcommand{\CA}{{\mathcal A}}
\newcommand{\CF}{{\mathcal F}}

\newcommand{\CU}{\mathcal{U}}
\newcommand{\CX}{\mathcal{X}}
\newcommand{\CC}{\mathcal{C}}

\newcommand{\I}{\mathcal{I}}

\newcommand{\CB}{\mathcal{B}}

\newcommand{\CL}{\mathcal{L}}
\newcommand{\CP}{\mathcal{P}}

\newcommand{\CD}{\mathcal{D}}
\newcommand{\CR}{\mathcal{R}}

\newcommand{\EE}{{\mathbb E}}
\newcommand{\PP}{{\mathbb P}}

\newcommand{\NN}{{\mathbb N}}
\newcommand{\rr}{{\mathbb R}}

\newcommand{\DD}{{\mathbb D}}

\newcommand{\FF}{{\mathbb F}}
\newcommand{\bS}{{\mathbb S}}

\newcommand{\VV}{{\mathbb V}}

\newcommand{\SL}{\mathscr{L}}

\newcommand{\SK}{\mathscr{K}}
\newcommand{\wdh}{\widehat}

\newcommand{\qv}[1]{\langle #1 \rangle}

\newcommand{\thmref}[1]{Theorem~{\rm \ref{#1}}}
\newcommand{\lemref}[1]{Lemma~{\rm \ref{#1}}}

\newcommand{\propref}[1]{Proposition~{\rm \ref{#1}}}
\newcommand{\defref}[1]{Definition~{\rm \ref{#1}}}
\newcommand{\remref}[1]{Remark~{\rm \ref{#1}}}

\newcommand{\assmref}[1]{Assumption~{\rm \ref{#1}}}
\newcommand{\secref}[1]{Section~{\rm \ref{#1}}}

\newcommand{\ito}{It\^{o} }

\newcommand{\gronwall}{Gr\"{o}nwall }

\title{Time inconsistent linear quadratic optimal control for mean field games in Hilbert spaces \thanks{This work was supported by a start-up fund from Wayne State University.}}
\author{Hongjiang Qian \thanks{Department of Mathematics, Wayne State University, Detroit, MI 48202 (\texttt{hjqian@wayne.edu}).}}
\date{}

\begin{document}

\maketitle

\begin{abstract}
This paper studies time-inconsistent linear-quadratic mean-field games governed by stochastic evolution equations in Hilbert spaces. The time inconsistency arises through the dependence of the cost functional on the initial time, while the mean-field interaction enters through the mean of the state distribution. By employing a multiperson differential game approach, we show that the closed-loop equilibrium strategy admits an affine feedback representation characterized by an operator-valued Riccati equation and a novel coupled forward-backward equilibrium system. We establish the existence and uniqueness of mild solutions to the equilibrium system and hence the existence of a closed-loop equilibrium strategy.
\end{abstract}



\section{Introduction}
Time-inconsistent stochastic control has attracted significant attention in recent years. In classical stochastic control, Bellman's principle of optimality provides the foundation for dynamic programming: a control that is optimal at the initial time remains optimal when the problem is restarted at a later time along the controlled trajectory. This property fails in many important models, for example, when the cost functional depends explicitly on the initial time, involves non-exponential discounting, or contains nonlinear dependence on conditional moments of the state. In such a situation, a precommitted optimal control is generally not dynamically consistent. The appropriate solution concept is therefore an equilibrium control, which is locally optimal against infinitesimal deviations and remains consistent with the decisions of future selves. We refer to Bj\"{o}rk-Khapko-Murgoci \cite{BKM17}, Hu-Jin-Zhou \cite{HJZ12, HJZ17}, Yong \cite{Yon12, Yon17}, and references therein for systematic accounts of time-inconsistent stochastic control.

Linear-quadratic (LQ) control constitutes one of the most tractable classes of stochastic control problems. In the time-consistent setting, the optimal feedback control is characterized by Riccati equations. The corresponding theory is well established both for stochastic differential equations in Euclidean spaces \cite{YZ99} and for stochastic evolution equations in Hilbert spaces \cite{GT14, Lu19}. For time-inconsistent LQ control, extensive developments have been made in finite-dimensional settings through equilibrium controls and equilibrium Riccati equations. In infinite-dimensional spaces, however, the theory remains much less developed. A notable contribution is the work of Dou and L\"{u} \cite{DL20}, who established closed-loop equilibrium strategies for time-inconsistent LQ control problems governed by stochastic evolution equations. The loss of compactness in infinite dimensions introduces substantial difficulties that are absent in finite-dimensional models, thus requiring refined a priori estimates.

Many practical systems involve interactions among a large number of agents whose decisions depend not only on their own states but also on the aggregate behavior of the population. Such problems naturally lead to the theory of mean-field games (MFGs), initiated independently by Lasry and Lions \cite{LL06a, LL06b, LL07} and Huang, Caines, and Malham\'e \cite{HCM07, HCM10}. MFGs describe the limiting behavior of Nash equilibria in large-population stochastic games through the optimization problem of a representative agent coupled with a consistency condition for the population distribution. In the time-consistent setting, this leads to a coupled system consisting of a Hamilton--Jacobi--Bellman (HJB) equation and a Fokker--Planck (or Fokker--Planck--Kolmogorov (FPK)) equation; see the monograph \cite{CD18}. When time inconsistency is present, the classical dynamic programming principle is no longer available, and the equilibrium structure becomes considerably more involved. For finite-dimensional time-inconsistent mean-field models, we refer to \cite{NZK18, Wan19, Yon17} for linear-quadratic problems and to \cite{MZ20, MNY25} for nonlinear mean-field control and mean-field games.

Recently, Federico, Ghilli, and Gozzi \cite{FGG25} studied time-consistent linear-quadratic mean-field games in Hilbert spaces. In their framework, the state is governed by a stochastic evolution equation driven by an additive cylindrical Wiener process, and the cost functional depends on both the individual state and the population mean. Unlike classical LQ control problems, the value function contains not only a quadratic component but also a first-order term and a constant term generated by the mean-field interaction. Consequently, the optimal feedback control has an affine structure and cannot be characterized solely by an operator-valued Riccati equation. Besides the Riccati equation for the quadratic part, one must analyze an additional forward-backward system describing the first-order component and the evolution of the population mean. Similar structures already arise in finite-dimensional mean-field control problems; see \cite[Chapter 6]{BFY13} and \cite{BSYY16}. We also refer to \cite{LF25, FGS25} for recent developments on mean-field games in Hilbert spaces.

Motivated by the developments in \cite{DL20} and \cite{FGG25}, a natural question is whether one can characterize closed-loop equilibrium strategies for time-inconsistent linear-quadratic mean-field games in Hilbert spaces. Despite the recent progress in both time-inconsistent control and mean-field games, this problem remains open in infinite-dimensional spaces. The purpose of the present paper is to address this question.

Specifically, we consider a class of time-inconsistent linear-quadratic MFGs governed by stochastic evolution equations in Hilbert spaces. The state equation is driven by a real-valued Brownian motion, while the cost functional is quadratic in the state, the control, and the mean of the state distribution. Time inconsistency arises through the dependence of the coefficients in the cost functional on the initial time. This framework encompasses a broad class of distributed-parameter systems, including models with spatial structure, age dependence, and path-dependent dynamics.

Our main contributions are threefold. First, we characterize closed-loop equilibrium strategies and show that they admit an affine feedback representation $u=-\Theta X+\Phi$.
Second, we derive the equations governing the feedback coefficients. The operator-valued component $\Theta$ is characterized by a backward Riccati-type equation, while $\Phi$ is characterized by a new coupled forward-backward equilibrium system. This additional equilibrium system reflects the interaction between time inconsistency and mean-field effects and cannot be recovered from the Riccati equation alone. Third, we establish the existence and uniqueness of mild solutions to the resulting equilibrium system, under a uniform bound on the Riccati equations of the multiperson game analogous to \cite[Assumption 3.2]{FGG25}, thereby proving the existence of a closed-loop equilibrium strategy.

Compared with the time-consistent Hilbert-space MFG theory in \cite{FGG25}, the present work introduces the equilibrium structure arising from time inconsistency. Compared with the time-inconsistent LQ theory for stochastic evolution equations in \cite{DL20}, it incorporates mean-field interactions and leads to a novel coupled forward-backward equilibrium system characterizing the affine feedback strategy. To the best of our knowledge, this is the first study of time-inconsistent linear-quadratic mean-field games governed by stochastic evolution equations in Hilbert spaces.

Our analysis is based on a multiperson differential game approximation. The time interval is partitioned into finitely many subintervals, and each player controls the system only on its own subinterval while taking the future players' equilibrium strategies as given. By passing to the limit as the mesh size tends to zero, we obtain the desired closed-loop equilibrium strategy. This construction naturally yields both the Riccati-type equation and the forward-backward equilibrium system that characterize the equilibrium feedback control.

The rest of the paper is organized as follows. \secref{sec:math} presents the mathematical formulation and the main results. \secref{sec:MFG} provides an overview of classical optimal control for mean-field games in Hilbert spaces. \secref{sec:mult-pers} develops the multiperson differential game approximation, while \secref{sec:conv} establishes its convergence. Building upon these results, \secref{sec:exist} derives the closed-loop equilibrium strategy.

\section{Mathematical Formulation and Main Results}\label{sec:math}
\subsection{Notation}
Let $(\Omega,\FF,\CF,\PP)$ be a complete filtered probability space, on which a one-dimensional standard Brownian motion $W$ is defined such that $\FF=\{\CF_t\}_{t\geq 0}$ is the natural filtration generated by $W$ augmented by all of $\PP$-null sets.  For any $[s,t]\subset [0,\infty)$ and a Banach space $\CX$, we let $L_{\CF_t}^2(\Omega;\CX)=\{\xi:\Omega \to \CX \, |\,  \xi \text{ is } \CF_t\text{-measurable}, \EE|\xi|_\CX^2 <\infty\}$ and
\bea
L_\FF^2([s,t];\CX)\ad\!\!\! =\{\phi:[s,t]\times \Omega \to \CX\,|\, \phi\cd \text{ is } \FF\text{-adapted}, \EE\int_s^t |\phi(s)|_\CX^2 ds <\infty\}, \\
C_\FF([s,t];L^2(\Omega;\CX))\ad\!\!\! =\{\phi:[s,t] \times \Omega \to \CX\,|\, \phi\cd\text{ is } \FF\text{-adapted,} \\
\aad\quad \text{ and } \phi:[s,t]\to L_{\CF_{t}}^2(\Omega;\CX) \text{ is continuous}\}.
\eea

Let $T>0$ and $t\in [0, T)$. Let $(H,\qv{\cdot,\cdot}_H)$ and $(U,\qv{\cdot,\cdot}_U)$ be two separable Hilbert spaces with norm $|\cdot|_H$ and $|\cdot|_U$ introduced by the corresponding inner products. By $\CP(H)$, we denote the space of regular probability measures on $H$. Moreover, we denote by $\CL(H)$ the space of linear bounded operators from $H$ to $H$, by $\bS(H)$ the subspace of $\CL(H)$ of self-adjoint operators, and by $\bS^+(H)$ the space of self-adjoint nonnegative operators. For any $0\leq t_1 \leq t_2 \leq T$, we denote by $C_s([t_1,t_2];\bS(H))$ the set of all strongly continuous mappings $\mathfrak{F}:[t_1,t_2] \mapsto \bS(H)$, that is, $\mathfrak{F}\cd\phi$ is continuous on $[t_1,t_2]$ for all $\phi\in H$, and $C_s([t_1,t_2];\CL(H))$ is defined in the same way. For brevity, we define $\SK_t:=L_{\CF_t}^2(\Omega;H)$. Moreover, we let
\bea
L_\FF^2(t_1,t_2;U)=\{\phi:[t_1,t_2] \times\Omega \to U \,|\, \phi \text{ is } \FF\text{-adapted, and } \EE\int_{t_1}^{t_2} |\phi(s)|_U^2 ds <\infty\}.
\eea
Throughout this paper, the norms on $H$, $U$, and the associated operator spaces are all denoted by $|\cdot|$. The distinction will be clear from the context.

\subsection{Mathematical
Formulation} Consider the following controlled stochastic evolution equations (SEEs): for any $s\in (t, T]$,
\beq{see}\left\{\barray\!\!\!
dX(s)= [(A+A_1(s))X(s) + B(s) u(s)]ds +[C(s)X(s)+ D(s) u(s)] dW(s), \\
\!\!\! X(t)=x \in L_{\CF_t}^2(\Omega;H).
\earray\right.\eeq
Here $A$ generates a $C_0$-semigroup on $H$, and $A_1\cd,C\cd\in L^\infty([0,T];\CL(H))$, and $B\cd,D\cd\in L^\infty([0,T]; \CL(U;H))$. For any $x\in \SK_t$ and $u\in \CU[t,T]:=L_\FF^2(t,T;U)$, the system \eqref{see} admits a unique solution $X\cd=X(\cdot\,;\,t,x,u\cd)\in C_\FF([t,T];L^2(\Omega;H))$; see \cite[Chapter 6]{DZ92}.
The cost functional we consider is
\beqa{cost}
J(t, x, u\cd) \ad = \EE_t \Big[ \int_t^T f(t,s,X(s),\pi(s),u(s)) ds + h(t,T, X(T), \pi(T)) \Big],
\eeqa
where $f:[0,T]\times [0,T] \times H\times \CP(H) \times U \to \rr$ and $h: [0,T] \times [0,T] \times H\times \CP(H) \to \rr$ are measurable functions given by
\bea
f(t,s,x,\pi,u)\ad := \half [\qv{Q(s,t) x, x}_H + \qv{R(s,t)u, u}_U \\
\aad \quad  + \qv{\bar Q(s,t) [x- \bar S(s,t) \int_H \xi\, \pi(s,d\xi)], x- \bar S(s,t) \int_H \xi\, \pi(s,d\xi)}],
\eea
\bea\!\!\!\disp
h(t,s,x,\pi)\!:=\!\half[\qv{F(t)x,x}_H \!+\qv{\bar F(t) [x-\bar G(t)\!\int_H\! \xi \pi(T,d\xi)],x - \bar G(t)\!\int_H\! \xi\pi(T,d\xi)}_H ],
\eea
where $Q,\bar Q,\bar S \in C([0,T]\times[0,T]; \bS(H)), R\in C([0,T]\times [0,T];\bS(U))$, and $F,\bar F, \bar G \in C([0,T];\bS(H))$ and for some $\dl>0$ and $0\leq t \leq s \leq T$,
\bea
Q(s,t)\wedge \bar Q(s,t)  \geq 0,\; R(s,t) \geq \dl I,\; F(t)\wedge \bar F(t) \geq 0.
\eea
In the sequel, the time variables $s$ and $t$ are suppressed in $A_1, B, C, D$ and $Q, R,\bar Q, \bar S$, $F, \bar F,\bar G$, if there is no confusion. To proceed, we make the following assumptions.
\begin{ass}\label{ass}
We assume
\begin{itemize}
\item[\rm{(H1)}] for any $0\leq t \leq \tau \leq s \leq T$, we have
\bea
\begin{aligned}\!\!\!\!
Q(s,t) & \leq Q(s,\tau), & \bar Q(s,t) & \leq \bar Q(s,\tau), &  \bar S(s,t) & \leq \bar S(s,\tau), &  R(s,t) & \leq R(s,\tau),\\
F(t) &\leq F(\tau), & \bar F(t) & \leq \bar F(\tau), &  \bar G(t) &\leq \bar G(\tau).
\end{aligned}
\eea
\item[\rm{(H2)}] $-\bar F(s) \bar G(s) \in \bS^+(H)$ and $-\bar Q(s,s)\bar S(s,s)\in \bS^+(H)$ for any $s\in [0,T]$ and the following conditions hold true: $\qv{\bar F(s) \bar G(s)x,x}_H=0 \Rightarrow \bar F(s)\bar G(s)x=0$ and $\qv{\bar Q(s,t)\bar S(s,t) x,x}_H=0 \Rightarrow \bar Q(s,t)\bar S(s,t) x=0$ for all $s,t\in [0,T]$ and $x\in H$.
\end{itemize}
\end{ass}
\begin{rem} (i)
A classical example that satisfies (H1) is the hyperbolic discount, e.g., $\zeta(s,t)=\zeta(s)/[a(s-t+1)]$ for $\zeta=Q,\bar Q,  R$, and $\ell(t)=\ell/[a(T-t+1)]$ for $\ell=F,\bar F$, where $a>0$, $F,\bar F$ are nonnegative, and $Q(s),\bar Q(s)\in \bS^+(H)$ and $R(s)\in \bS(U)$ such that $R(s)\geq \dl_1 I$ for some $\dl_1>0$; see \cite{DL20}. Moreover, $\bar S(s,t)\equiv-\sigma I$ and $\bar G(t)\equiv -\sigma' I$ with $\sigma,\sigma'\geq 0$ satisfy (H1), and together with the above $\bar Q$ and $\bar F$ they also satisfy (H2). (ii) Assumption (H2) is imposed to ensure the uniqueness of solutions to the forward-backward system \eqref{back-eta} in \secref{sec:MFG} over large time horizons; see \cite[Theorem 4.2]{FGG25}.
\end{rem}

The cost functional \eqref{cost} is well-defined. We study the following time-inconsistent stochastic linear quadratic (TISLQ) optimal control problem:
\smallskip

\noindent\textbf{Problem (TISLQ)}: find $\wdh u\in \CU[t,T]$ such that
\beq{tic}
J(t,x;\wdh u\cd)=\inf_{u\cd\in \CU[t,T]} J(t,x; u\cd)
\eeq

\begin{rem}
We note that if $Q(s,t)=Q(s), \bar Q(s,t)=\bar Q(s), \bar S(s,t)= \bar S(s)$, $R(s,t)=R(s)$, and $F(t)=F, \bar F(t)=\bar F, \bar G(t)= \bar G$ for any $0\leq t \leq s \leq T$, then the Problem (TISLQ) reduces to time-consistent LQ control for MFG in Hilbert spaces, which was recently investigated in \cite{FGG25} for SEEs driven by an additive cylindrical Wiener process. In contrast, the diffusion coefficient in the present work is state and control-dependent and therefore nonconstant.
\end{rem}

For any $\wdh u \in \CU[t, T]$  satisfying \eqref{tic}, it is called a \textit{precommitment} \textit{optimal control} for the Problem (TISLQ) at $(t,x)$. The corresponding $\wdh X$ and $(\wdh X,\wdh u)$ are called the \textit{precommitment optimal state} and \textit{precommitment optimal pair} of Problem (TISLQ), respectively. Due to time-inconsistency, for any given initial pair $(t,x)$, the precommitment optimal control $\wdh u$ will not be optimal afterward, thus it is not useful. In this paper, we instead focus on equilibrium controls, particularly closed-loop equilibrium strategies. Such strategies have been extensively studied in finite-dimensional settings; see, for example, \cite{HJZ12, HJZ17, Yon17}, and more recently in infinite-dimensional settings in \cite{DL20}. However, to the best of our knowledge, there are currently no results on time-inconsistent mean-field games in Hilbert spaces.

\subsection{Main Results}
Let $\Theta\in L^\infty([0,T];\CL(H;U))$ and $\Phi\in L^\infty([0,T];U)$. For any $(t,x)\in [0,T] \times \SK_t$, let $\wdh X\cd \equiv \wdh X(\cdot\,; t,x)$ satisfy the following equation:
\beq{sys-fed}
d\wdh X(s) =[(A+A_1-B \Theta)\wdh X(s) + B\Phi] ds + [(C-D\Theta) \wdh X(s)+D\Phi] dW(s), \forall\, s\in (t,T],
\eeq
with $\wdh X(t) =x \in \SK_t= L_{\CF_t}^2(\Omega;H)$
and for a.e. $s\in [t,T]$, define
\beq{u-hat}
\wdh u(s)=- \Theta(s) \wdh X(s) + \Phi(s).
\eeq
\begin{definition} The
$\Theta\cd$ and $\Phi\cd$ are called a closed-loop equilibrium strategy for Problem (TISLQ) if for any $t\in [0, T)$, $x\in \SK_t$ and $u\cd\in \CU[t, T]$,
\beq{eq-J}
\liminf_{\e\to 0^+} \frac{J(t,\wdh X(t); u^\e\cd)-J(t,\wdh X(t); \wdh u\cd)}{\e} \geq 0, \quad \PP\text{-a.s.}
\eeq
where $u^\e\cd=u\cd\indi_{[t,t+\e)}\cd + [-\Theta\cd X^\e\cd +\Phi\cd] \indi_{[t+\e,T]}\cd$, with $X^\e$ being the solution of \eqref{see} corresponding to the control $u^\e\cd$.
\end{definition}

To characterize the closed-loop equilibrium strategy, we introduce the following backward equation:
\beqa{bRic}\left\{\barray
-\partial_s \Ga(s,t)= \Ga(s,t)[A+A_1(s)-B(s)\Theta(s)] + [A+A_1(s)-B(s)\Theta(s)]^* \Ga(s,t) \\
\qquad\qquad\quad + [C(s)-D(s)\Theta(s)]^* \Ga(s,t) [C(s)-D(s)\Theta(s)]\\
\qquad\qquad\quad + [Q(s,t)+\bar Q(s,t)] +\Theta^*(s) R(s,t) \Theta(s), \quad 0\leq t\leq s \leq T, \\
\Ga(T,t) =F(t)+\bar F(t),\quad 0\leq t\leq T,
\earray\right.\eeqa
where $\partial_s$ denotes the derivative for $s$ variable and for a.e. $s\in [0,T]$, $\Theta$ is defined as
\beq{def-theta}
\Theta(s)= \big[R(s,s)+ D(s)^* \Ga(s,s) D(s) \big]^{-1}  \big[B(s)^* \Ga(s,s) + D(s)^* \Ga(s,s) C(s) \big].
\eeq
Moreover, we introduce the following forward-backward equations:
\beq{fb-eq}\left\{\barray
\ad\!\!\!\!\!\!\!\! \partial_s \ga(s,t)+[A+A_1(s)-B(s)\Theta(s)]^* \ga(s,t) +\Ga(s,t)B(s) \Phi(s) \\
\aad\!\!\!\!\!\!\!\!\quad\; + [C-D\Theta(s)]^* \Ga(s,t)D\Phi(s)-\Theta(s)^* R(s,t)\Phi(s)-\bar Q(s,t) \bar S(s,t) z(s) =0 \\
\ad \!\!\!\!\!\!\!\!\partial_s z(s) - [A+A_1(s)-B\Theta(s)] z(s) - B(s)\Phi(s) =0 \\
& \!\!\!\!\!\!\!\! \ga(T,t) =-\bar F(t)\bar G(t) z(T), \quad z(0)=z_0:=\int_H \xi\, \pi_0(d\xi)\in H,
\earray\right.\eeq
where for a.e. $s\in [0,T]$, $\Phi(s)$ is defined as
\beq{def-Phi}
\Phi(s)=-[R(s,s)+D(s)^* \Ga(s,s) D(s)]^{-1} B^{*} \ga(s,s).
\eeq

\begin{rem}
In \eqref{bRic} and \eqref{fb-eq}, $t$ is regarded as a parameter. That is, for any fixed $t\in [0,T]$, \eqref{bRic} and \eqref{fb-eq} are equations with the variable $s$. The equation \eqref{bRic} has been studied in \cite{DL20} except that we have extra $\bar Q$ and $\bar F$ in \eqref{bRic} coming from the mean-field interactions in the cost functional \eqref{cost}.  To our knowledge, the existence and uniqueness theory for \eqref{fb-eq}  is new in the literature and constitutes the main contribution of this paper. Moreover, this work lays the foundation for the study of time-inconsistent mean-field games with nonlinear dynamics and cost functionals and that of McKean--Vlasov-type stochastic partial differential equations.
\end{rem}

Let us define $
\Xi:= \big\{(s,t) \in [0,T] \times [0,T] \,|\, 0\leq t \leq s \leq T\big\}$, and $
\CX(\Xi):= \big\{\Ga(\cdot,\cdot) \in L^\infty(\Xi; \CL(H)) \,| \,  \Ga(\cdot,t) \in C_s([t,T]; \bS^+(H)) \text{ for a.e. }  t\in [0,T] \big\}$
\begin{definition}\label{def:mild}
Recall that $C_s([0, T];\bS^+(H))$ is the space of strongly continuous operator-valued functions.
\begin{itemize}
\item[\rm{(1)}] We say that $\Ga(\cdot,\cdot)\in \CX(\Xi)$ is a mild solution of \eqref{bRic} if for any $0\leq t \leq s \leq T$ and $\phi \in H$, we have
\beqa{mild-Ga}
\!\!\!\Ga(s,t)\phi\ad\!\!\!= \big\{ e^{(T-s)A^*} [F(t)+\bar F(t)] e^{(T-s)A} + \int_s^T e^{(\tau-s)A^*} \big\{ \Ga(\tau,t)(A_1-B\Theta) \\
\aad + (A_1-B\Theta)^* \Ga(\tau,t) + [Q(\tau,t)+\bar Q(\tau,t)] \\
\aad  + [C-D\Theta]^* \Ga(\tau,t) [C-D\Theta] +\Theta^*(\tau) R(\tau,t) \Theta(\tau)\big\} e^{(\tau-s)A} d\tau \big\}\phi.
\eeqa
\item[\rm{(2)}] Given $\Ga\in \CX(\Xi)$, we say that $\ga\in L^\infty(\Xi;H)$ with $\ga(\cdot,t)\in C([t,T];H)$ and $z\in C([0,T];H)$ are a mild solution of the forward-backward system \eqref{fb-eq}, if for any $0\leq t \leq s \leq T$,
\bea
\ga(s,t)\ad\!\!\!=e^{(T-s)A^*}(-\bar F(t)\bar G(t) z(T)) +\int_s^T e^{(\tau-s)A^*} \big[(A_1-B \Theta(\tau))^{*} \ga(\tau,t) \\
\aad\quad + \Ga(\tau,t) B\Phi(\tau)+(C-D\Theta(\tau))^{*}\Ga(\tau,t) D\Phi(\tau)\\
\aad\quad -\Theta(\tau)^* R(\tau,t)\Phi(\tau)-\bar Q(\tau,t)\bar S(\tau,t) z(\tau)\big]  d\tau, \\
z(s)\ad\!\!\! =e^{s A} z_0 + \int_0^s e^{(s-\tau)A}\big([A_1 -B\Theta(\tau)] z(\tau)+B\Phi(\tau)\big)d\tau.
\eea
\end{itemize}
\end{definition}

The main result of this paper is as follows.
\begin{theorem}\label{thm:main}
Suppose that \assmref{ass} and \assmref{ass-rho} hold. Then, there exist unique mild solutions $\Ga(\cdot,\cdot), (\ga(\cdot,\cdot),z(\cdot))$ of \eqref{bRic} and \eqref{fb-eq}, respectively, and a closed-loop equilibrium strategy $\Theta\cd$ and $\Phi\cd$ given by \eqref{def-theta} and \eqref{def-Phi}.
\end{theorem}
\assmref{ass-rho} is stated in \secref{sec:uni-bdd-PGa}. It is the analogue of \cite[Assumption 3.2]{FGG25} for the multiperson game. By \propref{prop:H3} and \remref{rem:H3}, it holds if $T$ is small, and for any $T>0$ if $\bar S$ and $\bar G$ are small.

\section{MFG in Hilbert Spaces}\label{sec:MFG}
In this section, we consider the time-consistent optimal control of MFG in Hilbert spaces. That is, we suppose that $Q(s,t)=Q(s), R(s,t)=R(s), \bar Q(s,t)= \bar Q(s), \bar S(s,t)= \bar S(s)$, $F(t)=F, \bar F(t)=\bar F$, and $\bar G(t)=\bar G$. Then, the running cost and terminal cost in \eqref{cost} become
\bea\!\!\!
f(x,\pi,u) \ad\!\!\!=\half\big[ \qv{(Q+\bar Q)x,x}_H -2\qv{\bar Q \bar S\int_H\xi\, \pi(\cdot,d\xi), x}_H + \qv{Ru,u}_U \\
\ad \quad + \qv{\bar Q \bar S \int_H\xi\, \pi(\cdot,d\xi), \bar S \int_H\xi\, \pi(\cdot,d\xi)}_H \big]\\
h(x,\pi) \ad \!\!\! = \half\big[ \qv{(F+\bar F)x, x}_H - 2\qv{\bar F \bar G \int_H\xi\, \pi(\cdot,d\xi), x}_H \\
\ad \quad + \qv{\bar F \bar G \int_H\xi\, \pi(\cdot,d\xi), \bar G \int_H\xi\, \pi(\cdot,d\xi)}_H \big].
\eea

Define the value function
\bea\ad
V(t,x)=\inf_{u\in \CU[t,T]}\EE_t\Big[\int_t^T f(X(s),\pi(s),u(s))ds + h(X(T),\pi(T))\Big].
\eea
Formally, if $V$ is smooth enough, we have the following HJB-FPK equation:
\beq{HJB}\left\{\barray
\disp \partial_t V(t,x)+\inf_{u\in \CU[t,T]} \big\{\CL^u V(t,x)+ f(x,\pi,u)\big\}=0,\\
\disp \partial_t \pi(t,dx)-\half\text{Tr}[(Cx+D\wdh{u})(Cx+D\wdh{u})^* \nabla^2 \pi(t,dx)] \\
\qquad\qquad\qquad +\text{div}(\pi(t,dx)((A+A_1)x+B\wdh u))=0 \\
V(T,x) =h(x,\pi(T)), \pi(0,dx)=\pi_0(dx),
\earray\right.\eeq
where the first equation is an infinite-dimensional parabolic PDE with $\CL^u$ given by
\bea\ad
\CL^u \psi(x): = \qv{(A+A_1)x+Bu, \nabla \psi(x)}_H +\half \text{Tr} \big[(Cx+Du)(Cx+Du)^* \nabla^2 \psi(x) \big],
\eea
for any $\psi:H\to \rr$. $\nabla$ and $\nabla^2$ denote the first and second order Fr\'{e}chet derivatives with respect to the $x$ variable in $H$. Given the path of $\pi\cd$, since $W$ is a real-valued Brownian motion, we have $
\text{Tr}[ (Cx+Du)(Cx+Du)^* \nabla^2 V(t,x)]=\qv{\nabla^2 V(t,x)(Cx+Du), Cx+Du}_H$.
Thus, the infimum over $u$ in \eqref{HJB} is attained at
\beq{feed-oc}
\wdh u= -[R+D^* \nabla^2 V(t, x) D]^{-1}[B^* \nabla V(t,x) + D^* \nabla^2 V(t,x)Cx],
\eeq
which is the control $\wdh u$ in the second equation of \eqref{HJB}.
Denoting by $\widehat{X}$ the optimally controlled state process, we impose the consistency condition $\SL_{\mathrm{law}}(\widehat{X}(s))=\pi(s)$ for all $s$. This condition can be reformulated as the Fokker-Planck-Kolmogorov equation for the distribution $\pi(\cdot)$, namely the second equation in \eqref{HJB}; see \cite{CD18}.

As this is an LQ control problem, we can guess $V(t,x)=\half[\qv{P(t)x,x}_H + 2\qv{\eta(t),x}_H + \theta(t)]$, for some $P:[0,T]\to \bS^+(H)$, $\eta:[0,T]\to H$, and $\theta:[0,T]\to \rr$. Then $\nabla V(t,x)= P(t)x+\eta(t), \nabla^2 V(t,x)= P(t)$. Substituting $\nabla V$ and $\nabla^2 V$ back to the optimal control \eqref{feed-oc} and the HJB equation \eqref{HJB}, the optimal control becomes
\beq{oc}
\wdh u =- K(t)^{-1}[L(t)x+B^*\eta(t)],
\eeq
where $
K(t):= R+D^* P(t) D$ and $L(t):= B^* P(t) + D^* P(t)C$, and we set $\Theta(t):=K(t)^{-1}L(t)$. Here, $P$ satisfies the following backward Riccati equation
\beq{P-Ric}\left\{\barray
-\dot{P} = P(A+A_1)+ (A+A_1)^* P + C^* P C + (Q+\bar Q) \\
\qquad -(PB + C^* P D) K^{-1} (B^* P + D^* P C),\qquad P(T)=F+\bar F,
\earray\right.
\eeq
and the $\eta$ is coupled with $z$ through the following forward-backward system:
\beq{back-eta}\left\{\barray
-\dot{\eta}(t) &\!\!\!= [A+A_1-B\Theta(t)]^* \eta - \bar Q \bar S z(t) \\
\dot{z}(t) &\!\!\!= [A+A_1 - B\Theta(t) ] z(t) - B K(t)^{-1} B^* \eta(t) \\
\eta(T) & \!\!\!= -\bar F\bar G z(T),\quad z(0) = \int_H \xi\, \pi(0,d\xi),
\earray\right.
\eeq
where $z\cd:=\int_H\xi\,\pi(\cdot,d\xi)$. For brevity, we omit the details for the derivation of \eqref{back-eta}, but refer to \cite{FGG25}. Moreover, the scalar function $\theta(t)$ has an explicit expression in terms of $P$ and $(\eta,z)$ as follows:
\bea
\theta(t)= \qv{\bar F \bar G z(T), \bar G z(T)}_H + \int_t^T \big[\qv{\bar Q\bar S z(s), \bar S z(s)}_H -  \qv{K^{-1}(s) B^* \eta(s), B^* \eta(s)}_U\big] ds.
\eea

The operator-valued Riccati equation \eqref{P-Ric} has been studied in \cite{Lu19, DL20} for LQ control problems of stochastic evolution equations. We have the following result.

\begin{lem}\label{lem:exi-P}
Under \assmref{ass}, there exists a unique mild solution $P$ satisfying \eqref{P-Ric} such that $0\leq P(s)\leq \Pi(s)$ for any $s\in [0,T]$, where $\Pi\in C_s([0,T];\bS^+(H))$ satisfies the following Lyapunov equation
\bea
- \dot{\Pi}(s)= \Pi(s) (A+A_1) + (A+A_1)^* \Pi(s) + C^* \Pi(s) C + (Q+\bar Q),
\eea
with the terminal value $\Pi(T)=F+\bar F$.
\end{lem}
\begin{proof}
The definition of a mild solution for $P$ is analogous to \defref{def:mild}. The proof follows from \cite[Theorem 2.2, Lemma 3.6]{Lu19}; see also \cite[Lemma 3.2]{DL20}. Consequently, the details are omitted.
\end{proof}

To establish the existence and uniqueness of solutions to the forward-backward system \eqref{back-eta}, we adopt the decoupling approach of \cite{FGG25} (cf. the finite-dimensional framework in \cite{BFY13}) and seek an ansatz of the form
\bea
\eta(t) = \rho(t)z(t),
\eea
for an operator-valued mapping $\rho: [0,T] \to \CL(H)$. Because the state equation \eqref{see} features a non-constant diffusion coefficient, \eqref{back-eta} differs structurally from the system in \cite{FGG25}. For completeness, we adapt the arguments of \cite{FGG25} to our setting. Both \eqref{P-Ric} and \eqref{back-eta} are analyzed via mild solutions; we omit the explicit definition here as it mirrors \defref{def:mild} \textit{mutatis mutandis} (see also \cite[Definition 2.2]{FGG25}).

\subsection{Existence of the solution for \eqref{back-eta}} Denote by $A^n \in \CL(H)$ the Yosida approximations of the operator $A$ as $A^n= n^2 \CR(n,A) - nI$, where $\CR(n,A)$ is the resolvent of  $A$. By \cite{BDDM07}, we have the following properties:
\beq{MT}
M_T := \sup_{t\in [0,T],n\in \NN}\|e^{tA^n}\|_{\CL(H)}=\sup_{t\in [0,T],n\in \NN}\|e^{t(A^n)^*}\|_{\CL(H)} <\infty
\eeq
and
\beqa{conv-An}\barray\begin{aligned}
e^{t A^n}x & \to e^{tA}x, & e^{t(A^n)^*}x & \to e^{tA^*}x,  \quad \forall\, t\in [0,T], x\in H \\
A^n x & \to A x,\ \forall\, x\in \CD(A),  & (A^n)^* x &\to A^* x, \quad \forall\, x\in  \CD(A^*).
\end{aligned}\earray
\eeqa
\begin{theorem}\label{thm:exist}
Suppose that \assmref{ass} holds. Then, the forward-backward system \eqref{back-eta} admits a mild solution.
\end{theorem}
\begin{proof}
The analysis proceeds by first decoupling the approximating system $(\eta_n, z_n)$ following \cite[Chapter 6]{BFY13}, and subsequently passing to the limit.
\smallskip

\noindent\textbf{Step 1}. Let us first consider the approximating system
\beq{back-etan}\left\{\barray
-\dot{\eta_n}(t) & \!\! = [A^n + A_1-B \Theta(t)]^* \eta_n(t) - \bar Q \bar S z_n(t) \\
\dot{z_n}(t) & \!\!= [A^n + A_1 - B\Theta(t)] z_n(t) - B K^{-1} B^* \eta_n(t) \\
\eta_n(T) & \!\!= -\bar F\bar G z_n(T),\quad z_n(0) = \int_H \xi\, \pi(0,d\xi),
\earray\right.
\eeq
We note that the boundedness of $A^n$ allows us to analyze the system via strict solutions. We postulate that a strict solution to \eqref{back-etan} exists in the form of
\beq{etan}
\eta_n(t)= \rho_n(t)z_n(t),
\eeq
where $\rho_n\in C_s^1([0,T];\bS(H))$, that is, $\rho_n\cd\phi\in C^1([0,T];H)$ for all $\phi\in H$, $\rho_n(T)=-\bar F\bar G$, and $z_n \in C^1([0,T];H)$. Differentiating \eqref{etan} yields $\dot{\eta}_n(t) = \dot{\rho_n}(t)z_n(t)+ \rho_n(t) \dot{z}_n(t)$.
Then, $\rho_n$ satisfies the following equation:
\beq{rho-n}
\dot{\rho_n}(t) =-[A^n+A_1-B\Theta]^*  \rho_n(t) - \rho_n(t) [A^n+A_1 - B \Theta] +\rho_n(t) B K^{-1}B^* \rho_n(t)+\bar Q\bar S,
\eeq
with $\rho_n(T)= - \bar F \bar G$,
while $z_n$ governs the evolution:
\beq{zn}
\dot{z_n}(t) = [(A^n+A_1)-B\Theta] z_n(t)- B K^{-1} B^* \rho_n(t)z_n(t),\quad
z_n(0) =z_0.
\eeq
Under \assmref{ass}, the Riccati equation \eqref{rho-n} admits a strict solution $\rho_n$, the proof of which parallels that of \cite[Proposition A.4]{FGG25}. Substituting this operator into \eqref{zn}, we readily obtain the unique strict solution $z_n$ to \eqref{zn}. Defining $\eta_n$ via the ansatz \eqref{etan}, a straightforward verification confirms that the pair $(\eta_n,z_n)$ indeed constitutes a strict solution to \eqref{back-etan}.
\smallskip

\noindent\textbf{Step 2}. Let $(\eta_n,z_n,\rho_n)$ be the functions defined in Step 1. To proceed, we provide uniform estimates for $(\eta_n,z_n)$, uniformly in $n$. As $z_n$ is a strict solution, it is also a mild solution. It then has the following integral form:
\beqa{int-zn}
z_n(t)\! = e^{tA^n} z_0 +\!\! \int_0^t e^{(t-s)A^n}(A_1- B \Theta) z_n(s)ds-\!\! \int_0^t  e^{(t-s)A^n} B K^{-1} B^* \rho_n(s)z_n(s)ds.
\eeqa

We first provide the estimate for $\Theta$.
Since $P(s)\geq 0$, we have $K(s)=R+D^*P(s)D\geq \dl I$, and thus $|K(s)^{-1}|_{\CL(U)}\leq \dl^{-1}$. Consequently, the definition of $\Theta$ implies that for any $s\in[0,T]$,
\beqa{bound-Theta}
|\Theta(s)|_{\CL(H;U)} \ad\!\!\! \leq |K(s)^{-1}|_{\CL(U)}\,|B^*P(s)+D^*P(s)C|_{\CL(H;U)} \\
\aad\!\!\! \leq \dl^{-1}\big(|B|_\infty+|D|_\infty|C|_\infty\big)|\Pi(s)|_{\CL(H)},
\eeqa
where the last inequality follows from \lemref{lem:exi-P}. Since $\Pi\in C_s([0,T];\bS^+(H))$, it is uniformly bounded, which consequently implies that $\Theta$ is also uniformly bounded. Given that $\rho_n$ is a strict solution of \eqref{rho-n}, (H2) gives $\rho_n(t)\geq0$ as in Step 1 of the proof of \cite[Proposition A.4]{FGG25}. So the term $\rho_nBK^{-1}B^*\rho_n\geq0$ can be dropped in the upper estimate of $\qv{\rho_n(t)x,x}_H$, and applying \eqref{MT} alongside the boundedness of $A_1$ and $B$ allows us to invoke Gronwall's inequality to establish the uniform boundedness of $\rho_n$. Thus, from \eqref{int-zn}, we obtain $
|z_n(t)|_H \leq c\, \big(|z(0)|_H + \int_0^t |z_n(s)|_H ds \big)$
which, by a second application of Gronwall's inequality, implies
\beq{est-zn}
|z_n(t)|_H \leq c, \quad \forall\,t\in [0,T].
\eeq
Define $\tilde z_n(t):=z_n(t)- e^{tA^n}z(0)$. It follows that $|\tilde z_n(t)|\leq c, \forall\, t\in [0,T]$. Moreover, from \eqref{int-zn}, for $0\leq s\leq t\leq T$ and $\psi\in H$, we have
\beqa{weak-equi}
|\qv{\tilde z_n(t)-\tilde z_n(s),\psi}_H|\ad=\Big|\int_s^t\qv{e^{(t-r)A^n}g_n(r),\psi}_Hdr+\qv{\tilde z_n(s),(e^{(t-s)(A^n)^*}-I)\psi}_H\Big|\\
\ad\leq c\,|t-s||\psi|_H+c\,|(e^{(t-s)(A^n)^*}-I)\psi|_H,
\eeqa
where $g_n:=(A_1-B\Theta)z_n-BK^{-1}B^*\rho_nz_n$.

We now provide estimates for $\eta_n$. From the mild representation of $\eta_n$, we have
\bea
\eta_n(t)& =e^{(T-t)(A^n)^*}(-\bar F \bar G z_n(T)) + \int_t^T e^{(s-t)(A^n)^*}[A_1 -B\Theta(s)]^* \eta_n(s) ds  \\
& \quad -\int_t^T e^{(s-t)(A^n)^*} \bar Q \bar S z_n(s) ds.
\eea
By \eqref{est-zn}, we have the estimate $
|\eta_n(t)|_H \leq c\, \big(|z_n(T)|_H + \int_t^T |\eta_n(s)|_H ds \big)$, which by the \gronwall inequality, implies that
\beq{est-etan}
|\eta_n(t)|_H \leq c, \quad \forall\, t\in [0,T].
\eeq
Furthermore, defining $\tilde \eta_n(t):= \eta_n(t)-e^{(T-t)(A^n)^*}(-\bar F\bar G z_n(T))$, we obtain $|\tilde \eta_n(t)|_H \leq c,\;\forall\, t\in [0,T]$ and in the same way $|\qv{\tilde\eta_n(t)-\tilde\eta_n(s),\psi}_H|\leq c\,|t-s||\psi|_H+c\,|(e^{|t-s|A^n}-I)\psi|_H$ for all $t,s\in[0,T]$ and $\psi\in H$.
\smallskip

\noindent\textbf{Step 3}. This step addresses the convergence of $(\eta_n, z_n)$ to a limit $(\eta, z)$. Let $\varrho > 0$ be a constant such that the sequences $\tilde{\eta}_n$ and $\tilde{z}_n$ from Step 2 take values in the closed ball $\mathcal{B}_\varrho :=\{x\in H: |x|_H \leq \varrho \}$,
for all $n\in \NN$. By \cite[Theorem 3.29]{Bre11}, the weak topology of a separable Hilbert space $H$ is metrizable on $\mathcal{B}_\varrho$. Let $\{\xi_k\}_{k=1}^\infty$ be a dense subset of $\mathcal{B}_\varrho$. We define the distance between $x$ and $y$ in $\mathcal{B}_\varrho$ by
\bea
d(x,y):=\sum_{k\in \NN}2^{-(k+1)} |\qv{x-y,\xi_{k}}_H|, \quad \forall\, x,y \in \CB_\varrho
\eea
This metric induces the weak topology on $\mathcal{B}_\varrho$ and satisfies the inequality $
d(x,y) \leq \varrho |x-y|_H$ for all $x,y \in \CB_\varrho.$

Since $(A^n)^*$ is the Yosida approximation of $A^*$, the convergences in \eqref{conv-An} are uniform in $t\in[0,T]$ \cite{BDDM07}. Since $|(e^{hA^n}-I)\psi|_H\leq\sup_{t\in[0,T]}|(e^{tA^n}-e^{tA})\psi|_H+|(e^{hA}-I)\psi|_H$ and each $A^n$ is bounded, we get $\sup_{n}|(e^{hA^n}-I)\psi|_H\to0$ as $h\to0^+$ for each $\psi\in H$, and likewise for $(A^n)^*$. Then \eqref{weak-equi}, its analogue for $\tilde\eta_n$ and the dominated convergence theorem for the series defining $d$ show that $\{\tilde z_n\}$ and $\{\tilde\eta_n\}$ are equicontinuous in $C([0,T];(\CB_\varrho,d))$.
Applying the Ascoli-Arzel\`{a} theorem in the metric space $C([0,T];(\CB_\varrho,d))$, we deduce the existence of a subsequence of $\{\tilde \eta_n, \tilde z_n\}_{n\in \NN}$ (still denoted by $n$ with certain abuse of notation) that converges to some $(\tilde \eta,\tilde z)\in C([0,T];(\CB_\varrho,d))$ such that
\bea\ad
\lim_{n\to \infty}\sup_{t\in [0,T]}\big[d(\tilde \eta_n, \tilde \eta)+ d(\tilde z_n, \tilde z)\big]=0.
\eea
Denote by $\rightharpoonup$ the weak convergence in $H$, this implies that $
\tilde \eta_n(t) \rightharpoonup \tilde \eta(t)$ and $ \tilde z_n(t) \rightharpoonup \tilde z(t)$ for all $t\in [0,T]$. For any $x\in H$, we have
\bea\ad
\lim_{n\to \infty}\sup_{t\in [0,T]} |e^{tA}x- e^{t A^n}x|_H =0.
\eea
This yields the weak convergence $z_n(t)= \tilde z_n(t) + e^{t A^n}z(0)  \rightharpoonup \tilde z(t)+ e^{tA}z(0)=: z(t)$ for all $t\in [0,T]$. Following a similar argument, we can establish the weak convergence of $\eta_n(t)$ to its corresponding limit $\eta(t)$. From the weak convergence of $z_n$, it follows that
\bea
e^{(T-t)(A^n)^*}(-\bar F \bar G z_n(T)) \rightharpoonup  e^{(T-t)A^*}(-\bar F \bar G z(T))
\eea
Then, the weak convergence of $\tilde \eta_n$ yields that for all $t\in [0,T]$
\bea \eta_n(t)= \tilde \eta_n(t) + e^{(T-t)(A^n)^*}[-\bar F \bar G z_n(T)] \rightharpoonup \tilde \eta(t)+ e^{(T-t)A^*}[-\bar F\bar G z(T)]=: \eta(t).
\eea

The preceding arguments establish the existence of a limit pair $(\eta, z)$. It remains to show that $(\eta, z)$ satisfies the equation \eqref{back-eta}. We omit the details for $\eta$, as they follow from an argument analogous to the one for $z$. Noting that $z_n$ is also a weak solution to \eqref{zn}, it follows that for any $\psi \in \mathcal{D}(A^*)$ and $t \in [0,T]$, we have
\beq{weak-zn}\barray
\qv{z_n(t),\psi}_H & =\qv{z(0),\psi}_H + \int_0^t\, \qv{z_n(s), (A^n)^* \psi}_H ds  \\
& + \int_0^t \, \qv{z_n(s), (A_1-B\Theta)^* \psi}_H ds- \int_0^t\, \qv{B K^{-1}B^* \eta_n(s), \psi}_H ds.
\earray\eeq
Taking the limit as $n \to \infty$ and using the weak convergence of $\eta_n$ and $z_n$, we obtain for any $t, s \in [0, T]$,
\bea
\qv{z_n(t),\psi}_H \ad \to \qv{z(t),\psi}_H \\
\qv{z_n(s), (A_1-B\Theta)^* \psi}_H \ad \to \qv{z(s), (A_1-B\Theta)^* \psi}_H \\
\qv{BK^{-1} B^* \eta_n(s),\psi}_H \ad \to \qv{B K^{-1} B^* \eta(s), \psi}_H.
\eea
Moreover, for each $s\in [0,T]$, as $n\to \infty$,
\bea\ad
|\qv{z_n(s), (A^n)^* \psi}_H -\qv{z(s), A^* \psi}_H|\\
\aad \leq |\qv{z_n(s), ((A^n)^* -A^*) \psi}_H|+ |\qv{z_n(s)-z(s), A^* \psi}_H| \to 0.
\eea
Indeed, by the weak convergence of $z_n$, we have $|\qv{z_n(s)-z(s), A^* \psi}_H |\to 0$. By uniform boundedness of $z_n$ in \eqref{est-zn} and \eqref{conv-An}, we have
\bea
|\qv{z_n(s), ((A^n)^* -A^*)\psi}_H| \leq |z_n(s)|_H |((A^n)^*-A^*)\psi|_H \leq c\,|((A^n)^* -A^*) \psi|_H \to 0.
\eea
Finally, we note that
\bea\ad
\sup_{s\in [0,T]}|\qv{z_n(s),(A^n)^* \psi}_H| \leq \sup_{s\in [0,T]}|z_n(s)|_H |(A^n)^* \psi|_H \leq c\,(1+|A^*\psi|_H).
\eea
The above upper bound is uniform in $n$. Thus, by the dominated convergence theorem, taking the limit as $n \to \infty$ in \eqref{weak-zn} yields
\bea
\qv{z(t),\psi}_H & =\qv{z(0),\psi}_H + \int_0^t \qv{z(s), A^* \psi}_H ds + \int_0^t \qv{z(s), (A_1-B\Theta)^* \psi}_H ds \\
& \quad -\int_0^t \qv{B K^{-1}B^* \eta(s), \psi}_H ds.
\eea
This indicates that $z$ is a weak solution. Since the equation for $z$ in \eqref{back-eta} is linear, $z$ is also a mild solution. This completes the proof of existence for \eqref{back-eta}.
\end{proof}

\subsection{Uniqueness of the solution for \eqref{back-eta}}
\begin{theorem}\label{thm:unique}
Let \assmref{ass} hold. Then the forward-backward system \eqref{back-eta} admits a unique solution.
\end{theorem}
\begin{proof}
Let $(\eta_1,z_1)$ and $(\eta_2,z_2)$ be two solutions to \eqref{back-eta}. Define $\wdh \eta:=\eta_1-\eta_2$ and $\wdh z:=z_1 -z_2$. It suffices to prove that $(\wdh \eta,\wdh z)=(0,0)$. Given $\wdh \eta$, the equation satisfied by $\wdh z$ is given by $
\dot{\wdh z}(t)= [A+A_1-B\Theta]\wdh z(t)- BK^{-1} B^* \wdh \eta(t), \wdh z(0)=0$,
and given $\wdh z$, the equation for $\wdh \eta$ is given by
\beq{eq-wdh-eta}
\dot{\wdh \eta}(t)+[A+A_1-B\Theta]^* \wdh \eta(t) -\bar Q \bar S \wdh z(t)=0, \quad \wdh \eta(T)=-\bar F \bar G \wdh z(T).
\eeq
By the continuity of $(\eta_j,z_j)$, $j=1,2$, there exists a constant $c> 0$ such that for all $t\in [0,T]$, $|\wdh \eta(t)|_H \vee |\wdh z(t)|_H \leq c$.

To our purpose, we consider $d\,\qv{\wdh \eta(t), \wdh z(t)}_H/dt$. As $A$ is unbounded, the \ito formula is not directly applicable. We consider the Yosida approximation $A^n$ of $A$. Let $(\wdh \eta_n ,\wdh z_n)$ satisfy the following approximating equation
\beq{wdh-etan}\left\{\barray
\dot{\wdh z}_n(t) \ad\!\!\!= [A^n+A_1-B\Theta] \wdh z_n(t)- B K^{-1} B^* \wdh \eta(t), \quad \wdh z_n(0)=0 \\
\dot{\wdh \eta}_n(t)\ad\!\!\! + [A^n+A_1-B\Theta]^* \wdh \eta_n(t)- \bar Q \bar S \wdh z(t)=0, \quad\wdh \eta_n(T)= -\bar F \bar G \wdh z(T).
\earray\right.\eeq

By \eqref{conv-An}, the dominated convergence theorem and Gronwall's inequality, one can readily  show that $\wdh \eta_n \to \wdh \eta$ and $\wdh z_n \to \wdh z$ uniformly on $[0,T]$ as $n\to \infty$. Omitting the details, this uniform convergence implies that $
\beta_n(t):=\qv{\wdh \eta_n(t), \wdh z_n(t)}_H \to \qv{\wdh \eta(t), \wdh z(t)}=:\beta(t)$, uniformly in $t\in [0,T]$.
By \eqref{wdh-etan}, we have
\bea
\dot{\beta}_n(t)=-\qv{B K^{-1}B^* \wdh \eta(t), \wdh \eta_n(t)}_H + \qv{\bar Q\bar S \wdh z(t), \wdh z_n(t)}_H.
\eea
From the convergence of $\wdh \eta_n$ and $\wdh z_n$, the right-hand side of above equation converges to
$
-\qv{BK^{-1}B^* \wdh \eta(t), \wdh \eta(t)}_H + \qv{\bar Q \bar S \wdh z(t), \wdh z(t)}_H$, uniformly on $[0,T]$. Therefore, we have
\beq{eq-eta-z}\barray
d\qv{\wdh \eta(t), \wdh z(t)}/dt \ad =  -\qv{K^{-1} B^* \wdh \eta(t), B^* \wdh \eta(t)}_{U} +\qv{\bar Q \bar S \wdh z(t), \wdh z(t)}_H.
\earray\eeq

By the assumption that $R\in \bS(U)$, $R\geq\dl I$ and $-\bar Q \bar S\in \bS^+(H)$, and noting that $\wdh z(0)=0$, we obtain
\bea\ad
d\qv{\wdh \eta(t), \wdh z(t)}_H/dt \leq 0, \quad  \qv{\wdh \eta(0), \wdh z(0)}_H=0.
\eea
Moreover, we have the terminal condition $
\qv{\wdh \eta(T), \wdh z(T)}_H = \qv{-\bar F \bar G \wdh z(T), \wdh z(T)}_H \geq 0$, which implies that $
\qv{\wdh \eta(t), \wdh z(t)}_H=0,\; \forall\, t\in [0,T]$. In particular, evaluating at $t=T$ yields $\qv{-\bar F\bar G \wdh z(T), \wdh z(T)}_H=0$. By our assumption on $\bar F\bar G$ in (H2), it follows that
\beq{ter-FG}
-\bar F\bar G \wdh z(T)=0.
\eeq
Moreover, from \eqref{eq-eta-z}, we have $\qv{\bar{Q} \bar S \wdh z(t), \wdh z(t)}_H=0$. The assumption on $\bar Q \bar S$ in (H2) then yields that
\beq{t-QS}
\bar Q \bar S \wdh z(t)=0, \quad \forall\,t\in [0,T].
\eeq
Substituting \eqref{ter-FG} and \eqref{t-QS} into \eqref{eq-wdh-eta} implies that $\wdh{\eta}(t) = 0$ for all $t \in [0, T]$, which in turn ensures that $\wdh{z}(t) = 0$. This completes the proof, establishing the uniqueness of the solution to the forward-backward system \eqref{back-eta}.
\end{proof}

\section{Multi-Person Differential Games}\label{sec:mult-pers}
We now apply the multiperson differential game approach to study the time inconsistency; see \cite{Yon17, DL20}. For $N\in \NN$, we introduce an associated $N$-person differential game. Define $
\Dl := \big\{t_k \,|\, 0\leq k \leq N\big\}\equiv \big\{ 0=t_0 < t_1 < t_2 <\cdots < t_{N-1} < t_N = T\big\}$ with its mesh size given by $
\|\Dl\| = \max_{0\leq k \leq N-1} (t_{k+1}-t_k)$. The player $k$ controls the system on  $[t_k, t_{k+1})$ with the cost functional being ``sophisticatedly'' constructed by the following
rules:
\begin{itemize}
    \item[(i)] each player plays optimally on the assumption that the later players play optimally;
    \item[(ii)] the $(k+1)$-th player's initial state is the $k$-th player's final state, and
    \item[(iii)] the $k$-th player will not be able to control the system from $t_{k+1}$ on, but she will ``discount'' the cost functional in her own way on the interval $[t_{k+1},T]$.
\end{itemize}
For $k=0,1,\dots,N-1$, we set
\bea
\begin{aligned}
Q_k(s) & = Q(s,t_k), &  \bar Q_k(s) & = \bar Q(s,t_k), &
 R_k(s) & =R(s,t_k), &  \bar S_k(s) & = \bar S(s,t_k) \\
 F_k & = F(t_k), & \bar F_k & = \bar F(t_k), & \bar G_k & = \bar G(t_k).
\end{aligned}
\eea
Moreover, for $t\in [t_k,t_{k+1}]$, we introduce the following notation
\beq{sub-k}
\begin{aligned}
K_k &:= R_k + D^* P_k D, & L_k &:= B^* P_k + D^* P_k C, & \Theta_k &:= K_k^{-1} L_k, & \Phi_k &:= -K_k^{-1} B^* \eta_k \\
\CA_k&:= A+ A_1- B \Theta_k, & \mathcal{C}_k & := C- D \Theta_k,
\end{aligned}
\eeq
where $P_k\in\CL(H),\eta_k \in H$ will be specified later.

\subsection{The Player $(N-1)$}
For player $(N-1)$, the player's state equation is
\beq{state-N-1}
\left\{\barray
dX_{N-1}(s) &\!\!\!\!\!\!\! = [(A+A_1)X_{N-1} + B u_{N-1}]ds \\
 &\!\!\!\!\!\!\! \quad +[C X_{N-1} + D u_{N-1}] dW(s),\, s\in [t_{N-1},t_N],  \\
X_{N-1}(t_{N-1}) &\!\!\!  = x_{N-1} \in \SK_{t_{N-1}},
\earray\right.
\eeq
and the cost functional is
\beqa{cost-N-1}\ad\!\!\!\!\!
J_{N-1}^\Dl(x_{N-1};u_{N-1}\cd) \\
\ad\!\!\!\!\! =\half\, \EE_{t_{N-1}}\Big\{
\int_{t_{N-1}}^{t_N} \Big[ \qv{Q_{N-1}X_{N-1}, X_{N-1}}_H + \qv{R_{N-1}u_{N-1}, u_{N-1}}_U  \\
\aad\!\!\!\!\! + \qv{\bar Q_{N-1} [X_{N-1}- \bar S_{N-1}\int_H \xi\, \pi(s,d\xi) ], X_{N-1}- \bar S_{N-1}\int_H \xi\, \pi(s,d\xi)} \Big] ds \\
\aad\!\!\!\!\! + \qv{F_{N-1} X_{N-1}(t_N), X_{N-1}(t_N)}_H + \qv{\bar F_{N-1} [X_{N-1}(t_N)- \bar G_{N-1} \int_H \xi\, \pi(t_N, d\xi)],\\
\aad\; X_{N-1}(t_N) - \bar G_{N-1} \int_H \xi\, \pi(t_N,d\xi)}
\Big\}.
\eeqa
Let us now consider the following optimal control problem:
\smallskip

\noindent{$\textbf{Problem (TISLQ)}_{N-1}$}. For $x_{N-1} \in \SK_{t_{N-1}}$, find a $\wdh u_{N-1}\cd\in \CU[t_{N-1}, t_N]$ such that
\beq{prob-N-1}
J_{N-1}^\Dl(x_{N-1}; \wdh u_{N-1}\cd)= \inf_{u_{N-1} \cd \in \, \CU[t_{N-1}, t_N]} J_{N-1}^\Dl(x_{N-1}; u_{N-1}\cd).
\eeq
This is a standard LQ problem for MFG in Hilbert spaces. Let $\wdh \pi_{N-1}=\SL_{\text{law}}(\wdh X_{N-1})$ and define $z_{N-1}\cd:=\int_H \xi\, \wdh \pi_{N-1}(\cdot,d\xi)$. By \secref{sec:MFG} and \eqref{sub-k}, there exists a $P_{N-1}\in C_s([t_{N-1},t_N];\bS^+(H))$ satisfying the following Riccati equation:
\beq{Ric-N-1}\left\{\barray
\dot{P}_{N-1}+ P_{N-1} (A+A_1)+ (A+A_1)^* P_{N-1}+C^* P_{N-1}C + [Q_{N-1}+\bar Q_{N-1}] \\
\qquad\;\;\, - (P_{N-1}B+ C^* P_{N-1}D ) K_{N-1}^{-1} (B^* P_{N-1}+ D^* P_{N-1}C ) = 0 \\
P_{N-1}(t_{N}) = F_{N-1}+\bar F_{N-1},
\earray\right.\eeq
and a $\eta_{N-1}$ coupled with $z_{N-1}$ through the following forward-backward system:
\beq{eta-N-1}\left\{\barray
- \dot{\eta}_{N-1} &\!\!\!\!= [A+A_1-B\Theta_{N-1} ]^* \eta_{N-1} -\bar Q_{N-1} \bar S_{N-1} z_{N-1}  \\
\dot{z}_{N-1} &\!\!\!\!= [A+A_1- B\Theta_{N-1}] z_{N-1} - B K_{N-1}^{-1} B^* \eta_{N-1},\\
\eta_{N-1}(t_N) &\!\!\!\! = - \bar F_{N-1} \bar G_{N-1} z_{N-1}(t_N),\; z_{N-1}(t_{N-1})= \int_H \xi\, \wdh \pi_{N-1}(t_{N-1},d\xi),
\earray\right.\eeq
where $\Theta_{N-1}$ and $K_{N-1}$ are defined in \eqref{sub-k}. The existence and uniqueness of solutions to \eqref{Ric-N-1} is analyzed in \secref{sec:MFG}. Following the approach in \secref{sec:MFG}, we can decouple $\eta_{N-1}$ as $\eta_{N-1}=\rho_{N-1}z_{N-1}$, where
\bea
\dot{\rho}_{N-1}=-(\CA_{N-1})^* \rho_{N-1}-\rho_{N-1}\CA_{N-1}+\rho_{N-1}B K_{N-1}^{-1} B^{*} \rho_{N-1}+\bar Q_{N-1}\bar S_{N-1},
\eea
with terminal condition $\rho_{N-1}(t_N)=-\bar F_{N-1}\bar G_{N-1}$. Since (H2) gives the sign of $\bar Q\bar S$ only on the diagonal, $-\bar Q_{N-1}\bar S_{N-1}$ need not be nonnegative. So the solvability of this equation is taken from (H3) in \secref{sec:uni-bdd-PGa} instead of \cite[Proposition A.4]{FGG25}. Given $\rho_{N-1}$, let $z_{N-1}$ solve $\dot z_{N-1}=[\CA_{N-1}-BK_{N-1}^{-1}B^*\rho_{N-1}]z_{N-1}$ with the initial value in \eqref{eta-N-1}, and set $\eta_{N-1}:=\rho_{N-1}z_{N-1}$. A direct computation shows that $(\eta_{N-1},z_{N-1})$ solves \eqref{eta-N-1}. For uniqueness, let $(\tilde\eta,\tilde z)$ be the difference of two solutions of \eqref{eta-N-1}. Arguing with the Yosida approximation as in \thmref{thm:unique}, $\mu:=\tilde\eta-\rho_{N-1}\tilde z$ satisfies $-\dot\mu=[\CA_{N-1}^*-\rho_{N-1}BK_{N-1}^{-1}B^*]\mu$ and $\mu(t_N)=0$. Hence $\mu=0$ and $\dot{\tilde z}=[\CA_{N-1}-BK_{N-1}^{-1}B^*\rho_{N-1}]\tilde z$ with $\tilde z(t_{N-1})=0$, which gives $\tilde z=0$ and $\tilde\eta=0$. Moreover, there exists a scalar-valued function $\theta_{N-1}$ satisfying
\beq{theta-1}\left\{\barray
\dot{\theta}_{N-1}+\qv{\bar Q_{N-1}\bar S_{N-1}z_{N-1}, \bar S_{N-1}z_{N-1}}_H -\qv{K_{N-1}^{-1}B^* \eta_{N-1}, B^* \eta_{N-1}}=0\\
\theta_{N-1}(t_N)=\qv{\bar F_{N-1}\bar G_{N-1}z_{N-1}(t_N),\bar G_{N-1}z_{N-1}(t_N)}_H.
\earray\right.
\eeq

Consequently, the optimal state $\wdh X_{N-1}\cd\equiv \wdh X_{N-1}(\cdot\,; t_{N-1}, x_{N-1})$ solves the closed-loop system:
\beq{wdhX-N-1}\left\{\barray
d\wdh X_{N-1}(s) = [(A+A_1-B\Theta_{N-1}) \wdh X_{N-1}(s) + B \Phi_{N-1} ] ds \\
\qquad\qquad\qquad + [(C-D\Theta_{N-1}) \wdh X_{N-1}(s) + D \Phi_{N-1} ]dW(s)\; \text{ in } [t_{N-1},t_N]\\
\wdh X_{N-1}(t_{N-1}) = x_{N-1},
\earray\right.
\eeq
and the optimal control admits the following feedback form:
\beq{oc-N-1}
\wdh u_{N-1}(s)= - \Theta_{N-1}(s) \wdh X_{N-1}(s) + \Phi_{N-1}(s), \;\text{ for a.e. } s\in [t_{N-1},t_N].
\eeq
Hence, player $(N-1)$ solved the Problem $\text{(TISLQ)}_{N-1}$.

\subsection{The Player $(N-2)$}
The state equation for the player $(N-2)$ is
\bea\left\{\barray
dX_{N-2} = [(A+A_1) X_{N-2}+  Bu_{N-2}]ds \\
\qquad\qquad + [CX_{N-2}+D u_{N-2} ] dW(s), \text{ in } [t_{N-2}, t_{N-1}]\\
X_{N-2}(t_{N-2})=x_{N-2} \in \SK_{t_{N-2}}.
\earray\right.
\eea
We emphasize that at $s = t_{N-1}$, player $(N-2)$ no longer controls the system. Instead, player $(N-1)$ takes over, applying her optimal control $\widehat{u}_{N-1}$ from \eqref{oc-N-1} on the interval $[t_{N-1}, t_{N}]$. Here, $\widehat{X}_{N-1}$ solves the closed-loop system \eqref{wdhX-N-1} with the initial state $x_{N-1} = X_{N-2}(t_{N-1})$. Consequently, the cost functional for player $(N-2)$ becomes
\beq{J-N-2}\barray\ad\!\!\!\!\!
J_{N-2}^\Dl (x_{N-2}; u_{N-2}\cd) \\
\aad\!\!\!\!\! =  \half\,\EE_{t_{N-2}} \Big\{ \int_{t_{N-2}}^{t_{N-1}} \big[ \qv{Q_{N-2}X_{N-2}, X_{N-2}}_H + \qv{R_{N-2} u_{N-2}, u_{N-2}}_U \\
\aad + \qv{\bar Q_{N-2}\, [X_{N-2}- \bar S_{N-2} \int_H\xi\, \pi(s,d\xi)], X_{N-2}-\bar S_{N-2} \int_H\xi\, \pi(s,d\xi)} \big] \, ds \\
\aad + \int_{t_{N-1}}^{t_N} \big[
\qv{Q_{N-2} \wdh X_{N-1}, \wdh X_{N-1}}_H + \qv{R_{N-2} \wdh u_{N-1}, \wdh u_{N-1}}_U \\
\aad  + \qv{\bar Q_{N-2}[\wdh X_{N-1}- \bar S_{N-2}\, z_{N-1}], \wdh X_{N-1}- \bar S_{N-2}\, z_{N-1}} \big] ds \\
\aad + \qv{F_{N-2}\,\wdh X_{N-1}(t_N), \wdh X_{N-1}(t_N)}_H  \\
\aad + \qv{\bar F_{N-2}\, [\wdh X_{N-1}(t_N)- \bar G_{N-2}\,  z_{N-1}(t_N)],\wdh X_{N-1}(t_N)- \bar G_{N-2}\, z_{N-1}(t_N)}
\Big\},
\earray\eeq
where we recall that $z_{N-1}\cd:=\int_H \xi\, \wdh \pi_{N-1}(\cdot,d\xi)$ where $\wdh \pi_{N-1}=\SL_{\text{law}}(\wdh X_{N-1})$.

We note that the discount factors in \eqref{J-N-2} over $[t_{N-1},t_N]$ are $Q_{N-2}, R_{N-2}, \bar Q_{N-2}$, $\bar S_{N-2}$ and $F_{N-2}, \bar F_{N-2}, \bar G_{N-2}$, respectively. This aligns with the definition of ``discounts in her own way'' for player $(N-2)$, consistent with the rules established at the beginning of this section.

To proceed, we need to rewrite the cost functional over $[t_{N-1}, t_N]$ in \eqref{J-N-2} such that player $(N-2)$ faces a standard linear-quadratic mean field game problem in Hilbert spaces on the interval $[t_{N-2}, t_{N-1}]$. To this end, for $t \in [t_{N-1}, t_N]$, we define
\bea\!\!
\VV_{N-2}(t,x)\ad\!\!\!\!\! := \half\,\EE_{t} \Big\{ \int_{t}^{t_N} \big[
\qv{Q_{N-2}\, \wdh X_{N-1}, \wdh X_{N-1}}_H + \qv{R_{N-2}\, \wdh u_{N-1}, \wdh u_{N-1}}_U \\
\aad\!\!\!\! + \qv{\bar Q_{N-2}[\wdh X_{N-1}- \bar S_{N-2}\,  z_{N-1}(s)], \wdh X_{N-1}- \bar S_{N-2}\, z_{N-1}(s) } \big] ds\\
\aad\!\!\!\! + \qv{F_{N-2}\, \wdh X_{N-1}(t_N), \wdh X_{N-1}(t_N)}_H \\
\aad\!\!\!\!  + \qv{\bar F_{N-2} [\wdh X_{N-1}(t_N)- \bar G_{N-2}\, z_{N-1}(t_N)], \wdh X_{N-1}(t_N)- \bar G_{N-2}\, z_{N-1}(t_N)}
\Big\},
\eea
with $\wdh X_{N-1}(t)=x$. By the definition of $\wdh u_{N-1}$ in \eqref{oc-N-1}, $\VV_{N-2}$ can be rewritten as
\bea\!\!\!
\VV_{N-2}(t,x) \ad\!\!\!\! = \half\, \EE_{t}\Big\{ \int_{t}^{t_{N}} \big[ \qv{(Q_{N-2}+ \Theta_{N-1}^* R_{N-2} \Theta_{N-1}+ \bar Q_{N-2})\wdh X_{N-1}, \wdh X_{N-1}}_H \\
\aad - 2 \qv{ \Theta_{N-1}^* R_{N-2} \Phi_{N-1} + \bar Q_{N-2} \bar S_{N-2}\,z_{N-1}, \wdh X_{N-1}}_H \\
\aad+ \qv{R_{N-2}  \Phi_{N-1},  \Phi_{N-1}}_U + \qv{\bar Q_{N-2} \bar S_{N-2}\, z_{N-1}, \bar S_{N-2}\, z_{N-1}}_H
\big] ds \\
\aad + \qv{(F_{N-2}+\bar F_{N-2})\, \wdh X_{N-1}(t_N), \wdh X_{N-1}(t_N)}_H  \\
\aad -2 \qv{\bar F_{N-2} \bar G_{N-2}\, z_{N-1}(t_N), \wdh X_{N-1}(t_N)}_H \\
\aad + \qv{\bar F_{N-2} \bar G_{N-2}\, z_{N-1}(t_N), \bar G_{N-2}\,z_{N-1}(t_N)}_H \Big\}.
\eea
Due to the LQ structure, we can guess that there exist deterministic functions
\bea
\Ga_{N-2}: [t_{N-1},t_N] \to \CL(H), \quad \ga_{N-2}: [t_{N-1},t_N] \to H, \quad
\chi_{N-2}: [t_{N-1},t_N] \to \rr
\eea
such that $
\VV_{N-2}(t,x)=\half[\qv{\Ga_{N-2}(t)x, x}_H +2 \qv{\ga_{N-2}(t),x}_H + \chi_{N-2}(t)]$.

Applying It\^{o}'s formula to $\VV_{N-2}(t,\wdh X_{N-1}(t))$, we find that $(\Ga_{N-2},\ga_{N-2}, \chi_{N-2}) $ satisfies the following equations on $[t_{N-1},t_{N}]$:
\beq{Ga-N-2}\left\{\barray &\!\!\!\!\!
\dot{\Ga}_{N-2}+ \Ga_{N-2} \CA_{N-1}+ \CA_{N-1}^* \Ga_{N-2} + \mathcal{C}_{N-1}^* \Ga_{N-2} \mathcal{C}_{N-1}  \\
&\quad\;\;\, + (Q_{N-2}+\bar Q_{N-2}) + \Theta_{N-1}^* R_{N-2} \Theta_{N-1}=0, \\
&\!\!\!\!\!\dot{\ga}_{N-2} + \CA_{N-1}^* \ga_{N-2} + \Ga_{N-2} B \Phi_{N-1} + \mathcal{C}_{N-1}^* \Ga_{N-2} D \Phi_{N-1}\\
&\quad\;\;\,  - \Theta_{N-1}^* R_{N-2} \Phi_{N-1} -\bar Q_{N-2} \bar S_{N-2}\, z_{N-1} =0, \\
&\!\!\!\!\! \dot{\chi}_{N-2}+ \qv{\Ga_{N-2} D\Phi_{N-1},D \Phi_{N-1}}_H + 2\qv{\ga_{N-2}, B \Phi_{N-1}}_H \\ &\quad\;\;\, + \qv{R_{N-2} \Phi_{N-1}, \Phi_{N-1}}_U + \qv{\bar Q_{N-2} \bar S_{N-2}\, z_{N-1}, \bar S_{N-2}  z_{N-1}}_H =0,
\earray\right.
\eeq
with conditions $
\Ga_{N-2}(t_N)= F_{N-2}+ \bar F_{N-2}$,
$
\ga_{N-2}(t_N) = -\bar F_{N-2} \bar G_{N-2}\, z_{N-1}(t_N) $,
and $\chi_{N-2}(t_N)= \qv{\bar F_{N-2} \bar G_{N-2}\, z_{N-1}(t_N), \bar G_{N-2} \, z_{N-1}(t_N)}_H$ at terminal time $t_N$. We observe that $\ga_{N-2}$ is coupled with $z_{N-1}$, the dynamics of which are specified in \eqref{eta-N-1}.

Therefore, the cost functional $J_{N-2}^\Dl(x_{N-2}; u_{N-2}\cd)$ in \eqref{J-N-2} can be written as
\bea\ad
J_{N-2}^\Dl(x_{N-2}; u_{N-2}\cd) \\
\aad
=\half\,\EE_{t_{N-2}} \Big\{ \int_{t_{N-2}}^{t_{N-1}}\big[ \qv{Q_{N-2}X_{N-2},X_{N-2}}_H + \qv{R_{N-2}\, u_{N-2}, u_{N-2}}_U \\
\aad\quad  + \qv{\bar Q_{N-2}(X_{N-2}- \bar S_{N-2}\int_H \xi\,\pi(s,d\xi)), X_{N-2} - \bar S_{N-2}\int_H\xi\,\pi(s,d\xi)}_H \big] ds \\
\aad\quad + \qv{\Ga_{N-2}(t_{N-1})\, X_{N-2}(t_{N-1}), X_{N-2}(t_{N-1})} \\
\aad\quad + 2\qv{\ga_{N-2}(t_{N-1}), X_{N-2}(t_{N-1})}_H + \chi_{N-2}(t_{N-1})\Big\}.
\eea
This expression matches the standard cost functional of an LQ MFG problem in Hilbert spaces on $[t_{N-2}, t_{N-1}]$. Player $(N-2)$ then faces the following control problem:
\smallskip

\noindent{\textbf{Problem (TISLQ)}$_{N-2}$}. For $x_{N-2} \in \SK_{t_{N-2}}$, find a $\wdh u_{N-2}\cd\in \CU[t_{N-2},t_{N-1}]$ s.t.
\beq{contr-N-2}
J_{N-2}^\Dl (x_{N-2}; \wdh u_{N-2} \cd) =\inf_{u_{N-2}\cd\in\, \CU[t_{N-2},t_{N-1}]} J_{N-2}^\Dl(x_{N-2}; u_{N-2}\cd).
\eeq
Similar to the case of player $(N-1)$, by virtue of \secref{sec:MFG}, there exists a mapping $P_{N-2}:[t_{N-2},t_{N-1}]\mapsto \CL(H)$ that satisfies the following Riccati equation:
\beq{P-N-2}\left\{\barray
\dot{P}_{N-2} + P_{N-2}(A+A_1) + (A+A_1)^* P_{N-2} + C^* P_{N-2}C + [Q_{ N-2}+\bar Q_{N-2}] \\
-(P_{N-2}B + C^* P_{N-2} D) (R_{N-2}+ D^* P_{N-2} D)^{-1} (B^* P_{N-2} + D^* P_{N-2}C) =0 \\
P_{N-2}(t_{N-1})= \Ga_{N-2}(t_{N-1}),
\earray\right.
\eeq
and a pair $(\eta_{N-2}, z_{N-2})$ that satisfies the following forward-backward equation:
\beq{eta-N-2}\left\{\barray
\dot{\eta}_{N-2} + [A+A_1-B\Theta_{N-2}]^*  \eta_{N-2} - \bar Q_{N-2}\bar S_{N-2}\,z_{N-2}=0 \;\text{ in } [t_{N-2},t_{N-1}] \\
\dot{z}_{N-2}=[A+A_1-B\Theta_{N-2}] z_{N-2} + B \Phi_{N-2}, \\
\eta_{N-2}(t_{N-1})=\ga_{N-2}(t_{N-1}),\quad z_{N-2}(t_{N-2})=\int_H \xi\, \wdh \pi_{N-2}(t_{N-2},d\xi),
\earray\right.\eeq
where $z_{N-2}(s):=\int_H\xi\,\wdh \pi_{N-2}(s,d\xi)$ and $\wdh\pi_{N-2}=\SL_{\text{law}}(\wdh X_{N-2})$. Proceeding as in \secref{sec:MFG}, we can decouple $\eta_{N-2}$ via the relation $\eta_{N-2}=\rho_{N-2}z_{N-2}$, where $\rho_{N-2}$ satisfies
\bea
\dot{\rho}_{N-2}=-(\CA_{N-2})^* \rho_{N-2} -\rho_{N-2}\CA_{N-2} +\rho_{N-2} B K_{N-2}^{-1}B^{*} \rho_{N-2}+\bar Q_{N-2}\bar S_{N-2},
\eea
with $\rho_{N-2}(t_{N-1})=\vartheta_{N-2}(t_{N-1})$, since $\eta_{N-2}(t_{N-1})=\ga_{N-2}(t_{N-1})$, where $\ga_{N-2}=\vartheta_{N-2}z_{N-1}$ on $[t_{N-1},t_N]$ as in \eqref{var-theta-k}.
The scalar-valued function $\theta_{N-2}$ satisfies
\beq{r-N-2}\left\{\barray\!\!
\dot{\theta}_{N-2} + \qv{\bar Q_{N-2} \bar S_{N-2}\,z_{N-2}, \bar S_{N-2}\, z_{N-2}}_H - \qv{K_{N-2}^{-1} B^* \eta_{N-2}, B^* \eta_{N-2}}=0\\
\!\! \theta_{N-2}(t_{N-1})= \chi_{N-2}(t_{N-1}).
\earray\right.
\eeq

Then, the optimal state process $\wdh X_{N-2}\cd$ solves the closed-loop system
\beq{X-N-2}\left\{\barray
d\wdh X_{N-2} = [ (A+A_1-B\Theta_{N-2})\wdh X_{N-2}  + B \Phi_{N-2} ]ds \\
\qquad\qquad  + [(C-D\Theta_{N-2}) \wdh X_{N-2} + D \Phi_{N-2} ] dW(s), \quad s\in [t_{N-2}, t_{N-1}) \\
\wdh X_{N-2}(t_{N-2})=x_{N-2},
\earray\right.
\eeq
and the optimal control takes the feedback form $\wdh u_{N-2}(s)= -\Theta_{N-2}(s) \wdh X_{N-2}(s) +\Phi_{N-2}(s)$.
Thus, the player $(N-2)$ solved Problem $\text{(TISLQ)}_{N-2}$. In contrast to player $(N-1)$, here we need to solve the three equations in \eqref{Ga-N-2} to reduce the complicated cost functional to a standard LQ-problem framework for MFGs in Hilbert spaces.

\subsection{The Player $(N-3)$} Following the same structure as the previous step, we now consider player $(N-3)$, for whom the state equation is
\beq{state-N-3}\left\{\barray
d X_{N-3}  = [(A+A_1) X_{N-3} + B u_{N-3}] ds \\
\qquad\qquad + [C X_{N-3}+ D u_{N-3}] dW(s),\quad s\in (t_{N-3},t_{N-2}], \\
X_{N-3}(t_{N-3}) = x_{N-3} \in \SK_{t_{N-3}}.
\earray\right.
\eeq

Because players $(N-2)$ and $(N-1)$ play optimally on $[t_{N-2},t_{N-1})$ and $[t_{N-1},t_N]$, respectively, the resulting optimal state process on $[t_{N-2},t_N]$, denoted by $X^\Dl$, satisfies
\bea\left\{\barray
dX^\Dl = \big[(A+A_1- B\Theta^\Dl) X^\Dl + B \Phi^\Dl \big] ds + \big[(C-D \Theta^\Dl) X^\Dl + D \Phi^\Dl \big] dW(s)\\
X^\Dl(t_{N-2}) = X_{N-3}(t_{N-2}),
\earray\right.
\eea
where for $s\in (t_{N-2},t_{N}]$, $\Theta^\Dl$ and $\Phi^\Dl$ are defined as
\bea
\Theta^\Dl(s) \ad =\Theta_{N-2}(s)\indi_{[t_{N-2},t_{N-1})}(s)+ \Theta_{N-1}(s) \indi_{[t_{N-1},t_N]}(s), \\
\Phi^\Dl(s) \ad = \Phi_{N-2}(s)\indi_{[t_{N-2},t_{N-1})}(s)+ \Phi_{N-1}(s) \indi_{[t_{N-1},t_N]}(s).
\eea
The optimal control on $[t_{N-2},t_N]$ is then given by $u^\Dl(s)=-\Theta^\Dl(s) X^\Dl(s)+ \Phi^\Dl(s)$. Moreover, we define the corresponding population mean of $X^\Dl(\cdot)$ as
\bea\ad
z^\Dl(s):=z_{N-2}(s)\indi_{[t_{N-2},t_{N-1})}(s)+ z_{N-1}(s) \indi_{[t_{N-1},t_N]}(s).
\eea
Consequently, $z^\Dl$ satisfies the following equation:
\beq{z-Dl}
\dot{z}^\Dl(s) = [A+A_1-B \Theta^\Dl] z^\Dl(s) + B \Phi^\Dl.
\eeq

The cost functional for the player $(N-3)$ is
\beq{cost-N-3}\barray\ad
J_{N-3}^\Dl (x_{N-3}; u_{N-3}\cd)\\
\aad =\half\, \EE_{t_{N-3}} \Big\{
\int_{t_{N-3}}^{t_{N-2}} \big[ \qv{Q_{N-3} X_{N-3}, X_{N-3}}_H + \qv{R_{N-3}\, u_{N-3}, u_{N-3}}_U \\
\aad\quad + \qv{\bar Q_{N-3} [X_{N-3}- \bar S_{N-3}\int_H \xi\, \pi(s,d\xi)], X_{N-3}-\bar S_{N-3}\int_H \xi\, \pi(s,d\xi)} \big] ds \\
\aad\quad +\int_{t_{N-2}}^{t_N} \big[\qv{Q_{N-3}\, X^\Dl, X^\Dl}_H + \qv{R_{N-3} \, u^\Dl, u^\Dl}_U \\
\aad\qquad \qquad \quad +\qv{\bar Q_{N-3}\,[X^\Dl(s) - \bar S_{N-3}\, z^\Dl(s)], X^\Dl(s)- \bar S_{N-3}\, z^\Dl(s)}_H \big] ds \\
\aad\quad + \qv{F_{N-3}\, X^\Dl(t_N), X^\Dl(t_N)}_H \\
\aad\quad + \qv{\bar F_{N-3}\, [X^\Dl(t_N)- \bar G_{N-3}\, z^\Dl(t_N)],  X^\Dl(t_N) - \bar G_{N-3}\, z^\Dl(t_N) }
\Big\}.
\earray\eeq
By a similar argument for player $(N-2)$, this cost functional rewrites as
\bea\ad
J_{N-3}^\Dl(x_{N-3}; u_{N-3}\cd)\\
\aad = \half\, \EE_{t_{N-3}} \Big\{\int_{t_{N-3}}^{t_{N-2}} \big[ \qv{Q_{N-3} X_{N-3}, X_{N-3}}_H+ \qv{R_{N-3}\, u_{N-3}, u_{N-3}}_U \\
\aad\quad+ \qv{\bar Q_{N-3}[X_{N-3} - \bar S_{N-3}\,\int_H \xi\, \pi(s,d\xi)], X_{N-3}- \bar S_{N-3}\,\int_H \xi\, \pi(s,d\xi)}_H \big] ds \\
\aad\quad + \qv{\Ga_{N-3}(t_{N-2}) X_{N-3}(t_{N-2}), X_{N-3}(t_{N-2})}_H \\
\aad\quad + 2\qv{\ga_{N-3}(t_{N-2}), X_{N-3}(t_{N-2})}_H + \chi_{N-3}(t_{N-2})
\Big\},
\eea
where the processes $\{\Ga_{N-3}, \ga_{N-3},\chi_{N-3}\}$ satisfy the following equations:
\beq{Ga-N-3}\left\{\barray
\dot{\Ga}_{N-3}+ \Ga_{N-3}(A+A_1 - B\Theta^\Dl) + (A+A_1-B\Theta^\Dl)^* \Ga_{N-3}+ (Q_{N-3}+\bar Q_{N-3}) \\
\qquad\;\; + (C-D\Theta^\Dl)^* \Ga_{N-3}(C-D\Theta^\Dl)
+ (\Theta^\Dl)^* R_{N-3} \Theta^\Dl =0,\quad s\in (t_{N-2},t_N],  \\
\Ga_{N-3}(t_N) = F_{N-3}+ \bar F_{N-3},
\earray\right.
\eeq
and
\beq{ga-N-3}\left\{\barray
\dot{\ga}_{N-3} + (A+A_1-B \Theta^\Dl)^{*} \ga_{N-3}+ \Ga_{N-3} B\Phi^\Dl + (C-D\Theta^\Dl)^* \Ga_{N-3} D \Phi^\Dl \\
\qquad\;\;- (\Theta^\Dl)^* R_{N-3} \Phi^\Dl -\bar Q_{N-3} \bar S_{N-3}\, z^\Dl =0,\qquad \qquad s \in (t_{N-2},t_N] \\
\dot{\chi}_{N-3} + \qv{\Ga_{N-3} D \Phi^\Dl,D\Phi^\Dl}_H + 2 \qv{\ga_{N-3}, B \Phi^\Dl}_H + \qv{R_{N-3} \Phi^\Dl,  \Phi^\Dl} \\
\qquad\;\; +\qv{\bar Q_{N-3} \bar S_{N-3}\, z^\Dl, \bar S_{N-3}\,z^\Dl}_H =0,
\earray
\right.
\eeq
with the terminal conditions $
\chi_{N-3}(t_N) = \qv{\bar F_{N-3} \bar G_{N-3} z^\Dl(t_N),  \bar G_{N-3} z^\Dl(t_N)}_H$ and  $\ga_{N-3}(t_N) = -\bar F_{N-3} \bar G_{N-3} z^\Dl(t_N)$. We note that $\ga_{N-3}$ is coupled with $z^\Dl$ in \eqref{z-Dl}.

Then, the optimal control problem for player $(N-3)$ is as follows:
\smallskip

\noindent{\textbf{Problem (TISLQ)}$_{N-3}$}. For $x_{N-3} \in \SK_{t_{N-3}}$, find a $\wdh u_{N-3}\cd\in \CU[t_{N-3},t_{N-2}]$ s.t.
\beq{J-N-3}
J_{N-3}^\Dl (x_{N-3}; \wdh u_{N-3} \cd) =\inf_{u_{N-3}\cd\in \CU[t_{N-3},t_{N-2}]} J_{N-3}^\Dl(x_{N-3}; u_{N-3}\cd).
\eeq
This is a standard LQ problem for MFG in Hilbert spaces. \secref{sec:MFG} implies that there exists a $P_{N-3}$ satisfying the following Riccati equation on $[t_{N-3},t_{N-2}]$:
\beq{P-N-3}\left\{\barray\!\!\!
\dot{P}_{N-3} + P_{N-3} (A+A_1) + (A+A_1)^* P_{N-3} + C^* P_{N-3}C + [Q_{ N-3}+\bar Q_{N-3} ] \\
\!\!\! -(P_{N-3}B + C^* P_{N-3} D) (R_{N-3} + D^* P_{N-3} D)^{-1} (B^* P_{N-3} + D^* P_{N-3}C) =0 \\
\!\!\! P_{N-3}(t_{N-2})= \Ga_{N-3}(t_{N-2}),
\earray\right.
\eeq
and a pair $(\eta_{N-3}, z_{N-3})$ satisfying the forward-backward equation:
\beq{eta-N-3}\left\{\barray
\dot{\eta}_{N-3} + [A+A_1-B\Theta_{N-3}]^* \eta_{N-3} - \bar Q_{N-3}\bar S_{N-3}\,z_{N-3}=0 \\
\dot{z}_{N-3} = [A+A_1 - B \Theta_{N-3} ] z_{N-3} + B \Phi_{N-3} \\
\eta_{N-3}(t_{N-2})=\ga_{N-3}(t_{N-2}),\quad z_{N-3}(t_{N-3}) = \int_H \xi\, \wdh \pi_{N-3}(t_{N-3},d\xi),
\earray\right.\eeq
where $z_{N-3}\cd:=\int_H\xi\,\wdh \pi_{N-3}(\cdot,d\xi)$ and $\wdh\pi_{N-3}=\SL_{\text{law}}(\wdh X_{N-3})$. The scalar function $\theta_{N-3}$ is then given by
\beq{r-N-3}\left\{\barray\!\!
\dot{\theta}_{N-3} + \qv{\bar Q_{N-3} \bar S_{N-3}\, z_{N-3}, \bar S_{N-3}\, z_{N-3}}_H - \qv{K_{N-3}^{-1} B^* \eta_{N-3}, B^* \eta_{N-3}}=0\\
\!\! \theta_{N-3}(t_{N-2})= \chi_{N-3}(t_{N-2}).
\earray\right.
\eeq
Similarly, we have $\eta_{N-3}=\rho_{N-3}z_{N-3}$, where $\rho_{N-3}$ satisfies
\bea
\dot{\rho}_{N-3} =-(\CA_{N-3})^* \rho_{N-3} -\rho_{N-3} \CA_{N-3}+ \rho_{N-3} B K_{N-3}^{-1} B^{*} \rho_{N-3}+\bar Q_{N-3}\bar S_{N-3},
\eea
with $\rho_{N-3}(t_{N-2})=\vartheta_{N-3}(t_{N-2})$.

Hence, the optimal state process $\wdh X_{N-3}\cd$ solves the closed-loop system
\beq{X-N-3}\left\{\barray
d\wdh X_{N-3} = [ (A+A_1-B\Theta_{N-3})\wdh X_{N-3}  + B \Phi_{N-3}]ds \\
\qquad\qquad + [(C-D\Theta_{N-3}) \wdh X_{N-3} + D \Phi_{N-3}] dW(s), \; s\in [t_{N-3}, t_{N-2}) \\
\wdh X_{N-3}(t_{N-3})=x_{N-3},
\earray\right.
\eeq
and the optimal control takes the feedback form $
\wdh u_{N-3}(s)= -\Theta_{N-3}(s) \wdh X_{N-3}(s) +\Phi_{N-3}(s)$. The player $(N-3)$ solved Problem $\text{(TISLQ)}_{N-3}$.

\subsection{Induction on Player $k$} Proceeding by induction, we construct the finite sequences $\{\Ga_k, \ga_k, \chi_k\}_{k=0}^{N-2}$ and $\{P_k, (\eta_k,z_k), \theta_k\}_{k=0}^{N-1}$. Specifically, $P_k$ is characterized by a Riccati equation, $(\eta_k,z_k)$ satisfies a coupled forward-backward system, and $\theta_k$ satisfies a scalar-valued equation. These equations take the following form:
\beq{P-k}\left\{\barray
\dot{P}_k + P_{k} (A+A_1) +(A+A_1)^* P_{k} + C^* P_{k}C + [Q_{ k}+\bar Q_{k} ] \\
\quad\;\, - (P_{k}B + C^* P_{k} D ) (R_{k}+ D^* P_{k} D)^{-1} (B^* P_{k} + D^* P_{k}C) =0 \; \text{ in } [t_k,t_{k+1}),\\
P_{k}(t_{k+1})= \Ga_{k}(t_{k+1}),
\earray\right.\eeq
with $\Ga_{N-1}(t_N)=F_{N-1}+\bar F_{N-1}$, and \beq{eta-k}\left\{\barray
\dot{\eta}_{k} + [A+A_1-B\Theta_k]^* \eta_{k}  - \bar Q_{k}\bar S_k\, z_k=0 \quad \text{ in } [t_k,t_{k+1}), \\
\dot{z}_k = [A+A_1-B\Theta_k] z_k + B \Phi_k \qquad\qquad\quad \text{in } [t_k,t_{k+1}),  \\
\eta_{k}(t_{k+1}) =\ga_{k}(t_{k+1}),\quad
z_k(t_k) = \int_H \xi\, \wdh \pi_k(t_k,d\xi),
\earray\right.\eeq
with $\ga_{N-1}(t_N)=-\bar F_{N-1}\bar G_{N-1}z^\Dl(t_N)$ and
\beq{r-k}\barray
\dot{\theta}_{k} + \qv{\bar Q_{k} \bar S_k\,z_k, \bar S_k\, z_k}_H - \qv{K_k^{-1} B^* \eta_{k}, B^*\eta_{k}}=0,\quad
\theta_{k}(t_{k+1})= \chi_{k}(t_{k+1})
\earray
\eeq
In the above, $z_k\cd:=\int_H\xi\,\wdh \pi_k(\cdot,d\xi)$ and $\wdh \pi_k=\SL_{\text{law}}(\wdh X_k)$, and $\chi_{N-1}(t_N)$ is given by the terminal condition below \eqref{Ga-ga-chi-k} with $k=N-1$.
As in \secref{sec:MFG}, we can decouple $\eta_k$ as $\eta_k=\rho_k z_k$, where $\rho_k$ satisfies
\beq{rho-k}
\dot{\rho}_k = -(\CA_k)^* \rho_k- \rho_k \CA_k + \rho_k B K_k^{-1} B^{*} \rho_k +\bar Q_k\bar S_k\quad \text{ in } [t_k,t_{k+1}),\quad \rho_k(t_{k+1})=\vartheta_k(t_{k+1}).
\eeq
As for player $(N-1)$, the solvability of \eqref{rho-k} is part of (H3) in \secref{sec:uni-bdd-PGa}, and $(\eta_k,z_k)$ is obtained from $\rho_k$ in the same way.

For any $s\in [0,t_N]$, we define
\beqa{Theta-Phi-Dl}\disp
\Theta^\Dl(s)=\sum_{k=0}^{N-1} \Theta_k(s) \indi_{(t_k,t_{k+1}]}(s), \; \Phi^\Dl(s)=\sum_{k=0}^{N-1} \Phi_k(s) \indi_{(t_k,t_{k+1}]}(s).
\eeqa
Similarly, we can define $\rho^\Dl, z^\Dl, K^\Dl, \CA^\Dl$, and $\mathcal{C}^\Dl$. At $s=0$, $\Theta^\Dl$, $\Phi^\Dl$ and these functions take the values of the functions with index $k=0$, and we use the convention $\vartheta_{N-1}(t_N):=-\bar F_{N-1}\bar G_{N-1}$.

Hence, the equations satisfied by $\rho^\Dl$ and $z^\Dl$ are given by
\beqa{rho-z-Dl}
\dot{\rho}^\Dl(s)\ad = -(\CA^\Dl)^* \rho^\Dl(s) - \rho^\Dl(s) \CA^\Dl + \rho^\Dl(s) B (K^\Dl)^{-1} B^{*} \rho^\Dl(s)+[\bar Q\bar S]^\Dl(s), \\
\dot{z}^\Dl(s) \ad = \CA^\Dl z^\Dl(s)- B(K^\Dl)^{-1} B^* \rho^\Dl(s) z^\Dl(s),
\eeqa
on each $(t_k,t_{k+1})$, with $\rho^\Dl(t_{k+1})=\vartheta_k(t_{k+1})$ for $0\leq k\leq N-1$, and $z^\Dl(0)=z_0$, where $[\bar Q\bar S]^\Dl(s):=\bar Q_k(s)\bar S_k(s)$ for $s\in(t_k,t_{k+1}]$. In general $\vartheta_k(t_{k+1})\neq\rho_{k+1}(t_{k+1})$, so $\rho^\Dl$ may jump at the nodes.

The sequence $\{\Ga_k,\ga_k, \chi_k\}_{k=0}^{N-2}$ defined on $[t_{k+1},t_N]$ satisfies the following equations:
\beq{Ga-ga-chi-k}\left\{\barray\!\!\!
\dot{\Ga}_k + \Ga_k \CA^\Dl + (\CA^\Dl)^* \Ga_k + (\mathcal{C}^\Dl)^* \Ga_k \mathcal{C}^\Dl + (Q_k+ \bar Q_k)+ (\Theta^\Dl)^* R_k \Theta^\Dl =0 \\
\!\!\!\dot{\ga}_k + (\CA^\Dl)^* \ga_k + \Ga_k  B \Phi^\Dl + (\mathcal{C}^\Dl)^* \Ga_k D \Phi^\Dl - (\Theta^\Dl)^* R_k  \Phi^\Dl -\bar Q_k \bar S_k z^\Dl = 0 \\
\!\!\! \dot{\chi}_k +\qv{\Ga_k D \Phi^\Dl, D\Phi^\Dl}\!+\! 2\qv{\ga_k, B\Phi^\Dl} \!+\! \qv{R_k \Phi^\Dl, \Phi^\Dl}_U\!+\! \qv{\bar Q_k \bar S_k\, z^\Dl, \bar S_k\,z^\Dl} =0,
\earray\right.\eeq
with terminal conditions:
\bea\!\!
\Ga_k(t_N) =F_k+\bar F_k,\,
\ga_k(t_N) =-\bar F_k \bar G_k z^\Dl(t_N),\,
\chi_k(t_N)  =\qv{\bar F_k \bar G_k z^\Dl(t_N), \bar G_k z^\Dl(t_N)}_H.
\eea

We consider the following closed-loop system
\beq{state-CD}\left\{\barray\!\!\!
dX^\Dl(s) \ad\!\!\! = \big[\CA^\Dl X^\Dl(s) + B \Phi^\Dl \big]ds +\big[\mathcal{C}^\Dl X^\Dl(s) + D \Phi^\Dl \big] dW(s), \;\text{ in } (0,t_N]\\
\!\!\! X^\Dl(0) \ad\!\!\! =x.
\earray\right.
\eeq
Let us define $
u^\Dl(s) = - \Theta^\Dl(s)X^\Dl(s)+ \Phi^\Dl(s), \text{ for a.e. } s\in [0,t_N]$.
Then for $0\leq k \leq N-1$, we have the control problem
\bea\ad
\inf_{u_k\cd\in \CU[t_k,t_{k+1}]} J_k^\Dl(X^\Dl(t_k); u_k\cd) = J_k^\Dl(X^\Dl(t_k); u^\Dl\cd|_{[t_k, t_N]}).
\eea

\section{Convergence Results}\label{sec:conv}
In this section, we establish the convergence of $\Theta^\Dl$ and $\Phi^\Dl$ defined in \eqref{Theta-Phi-Dl} as the mesh $\|\Dl\| \to 0$.

\subsection{Uniform Boundedness Estimates} \label{sec:uni-bdd-PGa}
Let $\tilde{Q},\tilde{F} \in \bS^+(H)$, and $\tilde{R} \in \bS^+(U)$ satisfy $
Q(s,s)+\bar Q(s,s)\leq  \tilde{Q},\; R(s,s)\leq \tilde{R},\; s\in [0,T],$ and $F(T)+ \bar F(T) \leq \tilde{F}$. By (H1), this gives $Q(s,t)+\bar Q(s,t)\leq \tilde Q$ and $R(s,t)\leq \tilde R$ for all $0\leq t\leq s\leq T$.

For any $0\leq k \leq N-1$, we consider the Riccati equation
\beq{wdh-P}\left\{\barray
\dot{\tilde{P}}_k + \tilde{P}_k (A+A_1)+ (A+A_1)^* \tilde{P}_k + C^* \tilde{P}_k C + \tilde{Q} \\
\quad\;\,- (\tilde P_k B+ C^* \tilde{P}_k D) (\tilde{R} + D^* \tilde{P}_k D)^{-1} (B^* \tilde{P}_k + D^* \tilde P_k C) =0 \;\text{ in }[t_k,t_{k+1}], \\
\tilde P_k(t_{k+1}) = \tilde \Ga_k(t_{k+1}),
\earray\right.
\eeq
and for $1\leq k \leq N-1$, consider the equation
\beq{wdt-Ga}\left\{\barray
\dot{\tilde{\Ga}}_{k-1} + \tilde\Ga_{k-1}(A+A_1-B\tilde \Theta_k) + (A+A_1-B\tilde \Theta_k)^* \tilde \Ga_{k-1} + \tilde Q\\
\qquad\;\, + (C- D \tilde \Theta_k)^* \tilde \Ga_{k-1} (C-D \tilde \Theta_k) + (\tilde \Theta_k)^* \tilde{R} \tilde \Theta_k =0\\
\tilde \Ga_{k-1}(t_N)= \tilde F,
\earray\right.
\eeq
where $
\tilde \Theta_k:= (\tilde{R} + D^* \tilde P_k D)^{-1} (B^* \tilde P_k + D^* \tilde P_k C)$  in $[t_k, t_{k+1}]$  for $0\leq k \leq N-1$,
and $\tilde\Ga_{N-1}(t_N):= \tilde F$. Moreover, we also consider the Lyapunov equation for $0\leq k\leq N-1$,
\beq{wdt-Pi}\left\{\barray
\dot{\tilde{\Pi}}_k + \tilde \Pi_k(A+A_1)+ (A+A_1)^* \tilde \Pi_k + C^* \tilde \Pi_k C + \tilde Q =0, \quad \text{ in } [t_k, t_{k+1}) \\
\tilde \Pi_k(t_{k+1})= \tilde \Ga_k(t_{k+1}).
\earray\right.
\eeq
Define the continuous time interpolation,
\bea\ad
\tilde\Pi^\Dl(s)= \sum_{k=0}^{N-1}\tilde\Pi_k(s) \indi_{(t_k,t_{k+1}]}(s), \quad s\in [0,T].
\eea
Then by \cite[Proposition 4.1]{DL20}, we have the following comparison results.
\begin{lem}\label{lem:compare}
Let assumption (H1) hold. Then we have $
0\leq \Ga_\varsigma (s) \leq P_k(s) \leq \tilde P_k(s) \leq \tilde \Pi_k(s)$ for all $0\leq \varsigma <k \leq N-1$,
with the convention that $
\Ga_{-1}(s)=0$. Moreover, $\tilde \Pi^\Dl\cd$ is uniformly bounded for all partitions $\Dl$ on $[0,T]$.
\end{lem}
\begin{proof}
The proof of this lemma can be found in \cite[Propositions 4.1 and 4.2, pp. 498--502]{DL20}, details are omitted.
\end{proof}

The preceding lemma ensures that $P_k(\cdot)$ and $\Ga_k(\cdot)$ are uniformly bounded with respect to both $k$ and the partition $\Dl$. To establish the uniform boundedness of $\{\eta_k, \theta_k\}$, by an argument analogous to \eqref{bound-Theta} and leveraging the uniform bound on $P_k$, we obtain the uniform boundedness of $\Theta_k$. Together with the bound of $\rho_k$ in (H3) below, this gives the uniform boundedness of $z_k$ and then that of $\eta_k=\rho_kz_k$. The remaining details are omitted for brevity.

We now establish the uniform boundedness of $\{\ga_k,\chi_k\}$. In the spirit of \secref{sec:MFG}, we handle $\ga_k$ by introducing the decoupling ansatz $\ga_k(s)=\vartheta_k(s)z^\Dl(s)$, where $\vartheta_k:[t_{k+1},T] \to \CL(H)$. Differentiating this relation with respect to time yields $\dot{\ga}_k=\dot{\vartheta}_k z^\Dl + \vartheta_k \dot{z}^\Dl$, which allows us to derive the evolution equation for $\vartheta_k$ as follows:
\beq{var-theta-k}
\dot{\vartheta}_k+(\CA^\Dl)^* \vartheta_k + \vartheta_k \CA^\Dl + \vartheta_k B \Psi^\Dl+
[ \Ga_k B + (\CC^\Dl)^* \Ga_k D- (\Theta^\Dl)^* R_k]\Psi^\Dl -\bar Q_k \bar S_k =0,
\eeq
with terminal condition $\vartheta_k(t_N)=-\bar F_k \bar G_k$, where $\Psi^\Dl :=-(K^{\Dl})^{-1}B^* \rho^\Dl$. Since $\Phi^\Dl=\Psi^\Dl z^\Dl$, the Yosida approximation as in \thmref{thm:unique} shows that $\vartheta_kz^\Dl$ is a mild solution of the linear equation for $\ga_k$ in \eqref{Ga-ga-chi-k}, and hence $\ga_k=\vartheta_kz^\Dl$ by uniqueness. It has the following mild representation:
\beqa{mild-vartheta}\disp
\vartheta_k(s) =e^{(T-s)A^*}(- \bar F_k \bar G_k) e^{(T-s)A} + \int_s^T e^{(r-s)A^*}\{(A_1-B\Theta^\Dl)^* \vartheta_k(r) \\
\disp +\vartheta_k(r)(A_1-B\Theta^\Dl) + \vartheta_k(r) B \Psi^\Dl + [\Ga_k B+ (\CC^\Dl)^* \Ga_k D- (\Theta^\Dl)^* R_k]\Psi^\Dl - \bar Q_k \bar S_k\}e^{(r-s)A}dr.
\eeqa
By \lemref{lem:compare} and the uniform boundedness of $\Theta^\Dl,\Psi^\Dl$, and $z^\Dl$, there exists a constant $c>0$ such that $\sup_{s\in [t_{k+1},T]}|\vartheta_k(s)|_{\CL(H)} \leq c$, thus $|\ga_k|_{L^\infty([0,T];H)} \leq c$.  Analogously, we can prove the uniform boundedness of $\chi_k$ by using that of $\Ga_k$ and $\ga_k$.

We now turn to $\rho_k$. Unlike \secref{sec:MFG}, the terminal value $\rho_k(t_{k+1})=\vartheta_k(t_{k+1})$ need not belong to $\bS^+(H)$, since \eqref{var-theta-k} contains the nonsymmetric terms $\vartheta_kB\Psi^\Dl$ and $\La_k\Psi^\Dl$ with
\beq{La-k}
\La_k:=\Ga_kB+(\CC^\Dl)^*\Ga_kD-(\Theta^\Dl)^*R_k,
\eeq
and (H2) gives the sign of $\bar Q\bar S$ only on the diagonal. So \cite[Proposition A.4]{FGG25} does not apply to $\rho_k$. As in \cite[Assumption 3.2]{FGG25}, we impose the following condition.
\begin{ass}\label{ass-rho}
\begin{itemize}
\item[\rm{(H3)}] There exists $c_\rho>0$ such that for any partition $\Dl$ of $[0,T]$ and $0\leq k\leq N-1$, \eqref{rho-k} admits a mild solution $\rho_k\in C_s([t_k,t_{k+1}];\CL(H))$ with $|\rho_k(s)|_{\CL(H)}\leq c_\rho$ for all $s\in[t_k,t_{k+1}]$.
\end{itemize}
\end{ass}
Let $c_\Theta:=\sup_\Dl|\Theta^\Dl|_\infty$ and $c_\La:=\sup_\Dl\max_{k}\sup_{s\in[t_{k+1},T]}|\La_k(s)|$, which are finite by \lemref{lem:compare} and \eqref{bound-Theta}. Set
\beq{c012}
c_0:=|\bar F\bar G|_\infty+T|\bar Q\bar S|_\infty,\quad c_1:=2\big(|A_1|_\infty+|B|_\infty c_\Theta\big)+\dl^{-1}|B|_\infty c_\La,\quad c_2:=\dl^{-1}|B|_\infty^2,
\eeq
where $|\bar F\bar G|_\infty:=\sup_{t\in[0,T]}|\bar F(t)\bar G(t)|_{\CL(H)}$ and $|\bar Q\bar S|_\infty:=\sup_{(s,t)\in\Xi}|\bar Q(s,t)\bar S(s,t)|_{\CL(H)}$.

\begin{prop}\label{prop:H3}
Let \assmref{ass} hold and assume that
\beq{H3-cond}
M_T^4c_0c_2\int_0^T e^{M_T^2c_1r}dr<1.
\eeq
Then (H3) holds with $c_\rho:=M_T^2c_0e^{M_T^2c_1T}\big[1-M_T^4c_0c_2\int_0^Te^{M_T^2c_1r}dr\big]^{-1}$.
\end{prop}
\begin{proof}
Fix $\Dl$ and $k$, and suppose that $\rho_{k+1},\dots,\rho_{N-1}$ have been constructed. Then $\Psi^\Dl$ is bounded on $[t_{k+1},T]$ and $\vartheta_k$ is given by the linear equation \eqref{mild-vartheta}. Let $\rho_k$ be a mild solution of \eqref{rho-k} on $[\sigma,t_{k+1}]$ for some $\sigma\in[t_k,t_{k+1}]$, and set $\wdh\vartheta_k:=\rho_k$ on $[\sigma,t_{k+1}]$ and $\wdh\vartheta_k:=\vartheta_k$ on $(t_{k+1},T]$. Since $\Psi^\Dl=-K_k^{-1}B^*\rho_k$ on $(t_k,t_{k+1}]$, we have $\rho_kB\Psi^\Dl=-\rho_kBK_k^{-1}B^*\rho_k$ there. By the semigroup property, the mild forms of $\rho_k$ and $\vartheta_k$ combine into
\beqa{concat}
\wdh\vartheta_k(s)\ad=e^{(T-s)A^*}(-\bar F_k\bar G_k)e^{(T-s)A}+\int_s^Te^{(r-s)A^*}\big\{(A_1-B\Theta^\Dl)^*\wdh\vartheta_k(r)+\wdh\vartheta_k(r)(A_1-B\Theta^\Dl)\\
\aad+\wdh\vartheta_k(r)B\Psi^\Dl+\indi_{(t_{k+1},T]}(r)\La_k\Psi^\Dl-\bar Q_k\bar S_k\big\}e^{(r-s)A}dr,\quad s\in[\sigma,T].
\eeqa
The same holds for $\wdh\vartheta_j$ on $[t_j,T]$, $k<j\leq N-1$, with $\wdh\vartheta_{N-1}=\rho_{N-1}$. Define
\[
\upsilon(s):=\max_{k\leq j\leq N-1}\ \sup_{r\in[s\vee t_j,T]}|\wdh\vartheta_j(r)|_{\CL(H)},\quad s\in[\sigma,T].
\]
Then $\upsilon$ is finite and nonincreasing. For $r\in(t_j,t_{j+1}]$ we have $\rho^\Dl(r)=\wdh\vartheta_j(r)$, hence $|\Psi^\Dl(r)|\leq\dl^{-1}|B|_\infty\upsilon(r)$. Thus \eqref{concat} and \eqref{MT}, which also gives $|e^{tA}|_{\CL(H)}\leq M_T$ by \eqref{conv-An}, yield
\beq{ineq-ups}
\upsilon(s)\leq M_T^2c_0+M_T^2\int_s^T\big[c_1\upsilon(r)+c_2\upsilon(r)^2\big]dr=:Y(s),\quad s\in[\sigma,T].
\eeq
We may assume $c_0>0$, so that $Y\geq M_T^2c_0>0$. With $\lambda:=M_T^2c_1$, we get for a.e. $s$
\[
\frac{d}{ds}\frac{e^{-\lambda s}}{Y(s)}=e^{-\lambda s}\Big[\frac{M_T^2(c_1\upsilon(s)+c_2\upsilon(s)^2)}{Y(s)^2}-\frac{\lambda}{Y(s)}\Big]\leq M_T^2c_2e^{-\lambda s}.
\]
Integrating over $[s,T]$, we obtain
\[
\frac{1}{Y(s)}\geq \frac{e^{-\lambda(T-s)}}{M_T^2c_0}-M_T^2c_2\int_0^{T-s}e^{-\lambda r}dr\geq \frac{e^{-\lambda T}}{M_T^2c_0}-M_T^2c_2\int_0^{T}e^{-\lambda r}dr=\frac{1}{c_\rho},
\]
and $c_\rho>0$ by \eqref{H3-cond}. Hence $\upsilon\leq Y\leq c_\rho$ on $[\sigma,T]$. In particular, every mild solution of \eqref{rho-k} on $[\sigma,t_{k+1}]$ satisfies $|\rho_k|\leq c_\rho$. With this a priori bound, the contraction of \cite[Proposition A.2]{FGG25} on balls of radius $2M_T^2c_\rho$ has a step length depending only on $c_\rho$, and iterating it as in Step 2 of the proof of \cite[Proposition A.4]{FGG25} gives the mild solution $\rho_k$ on $[t_k,t_{k+1}]$ with $|\rho_k|\leq c_\rho$. Backward induction on $k$ yields (H3).
\end{proof}

\begin{rem}\label{rem:H3}
Condition \eqref{H3-cond} holds if $T$ is small, which corresponds to \cite[Proposition A.2]{FGG25}. The constants $c_1$ and $c_2$ do not depend on $\bar S$ and $\bar G$, while $c_0\to0$ as $\bar S,\bar G\to0$. So \eqref{H3-cond} also holds for any $T>0$ if $\bar S$ and $\bar G$ are small, that is, if the mean-field coupling is weak. In the time-consistent case of \secref{sec:MFG}, we have $\Ga_k=P^\Dl$ on $(t_{k+1},T]$, and thus $\La_k=(L^\Dl)^*-(\Theta^\Dl)^*K^\Dl=0$ with $L^\Dl:=B^*P^\Dl+D^*P^\Dl C$. Then $\vartheta_k$ and $\rho^\Dl$ coincide with $\rho$ in \secref{sec:MFG}, and (H3) holds for any $T>0$ by (H2) and \cite[Proposition A.4]{FGG25}.
\end{rem}

\subsection{Convergence Analysis}\label{sec:conv-PGa-gachi}
Given a partition $\Dl$ of $[0,T]$, for $s\in [0,T]$, define
\bea
Q^\Dl(s) \ad =\sum_{k=0}^{N-1} Q(s,t_k) \indi_{(t_k,t_{k+1}]}(s), \quad R^{\Dl}(s)= \sum_{k=0}^{N-1} R(s,t_k) \indi_{(t_k,t_{k+1}]}(s), \\
\bar Q^\Dl(s) \ad = \sum_{k=0}^{N-1} \bar Q(s,t_k) \indi_{(t_k,t_{k+1}]}(s), \quad \bar S^\Dl(s)= \sum_{k=0}^{N-1} \bar S(s,t_k) \indi_{(t_k,t_{k+1}]}(s).
\eea
For $\ell=P,\eta$ and $\theta$, we define
\beq{def-ell-Dl}
\ell^\Dl(s)=  \sum_{k=0}^{N-1}\ell_k(s) \indi_{(t_k,t_{k+1}]}(s),\quad  s \in (0,T],\quad \ell^\Dl(0):=\ell_0(0).
\eeq
For $0\leq k\leq N-1$, let
\beq{hat-def}
(\wdh\Ga_k,\wdh\ga_k,\wdh\chi_k,\wdh\vartheta_k):=\begin{cases}(P_k,\eta_k,\theta_k,\rho_k)&\text{on }[t_k,t_{k+1}],\\ (\Ga_k,\ga_k,\chi_k,\vartheta_k)&\text{on }(t_{k+1},T].\end{cases}
\eeq
For $\ell=\Ga,\ga$, $\chi$ and $\vartheta$, we define
\beq{def-ell-Dl-s}
\ell^\Dl(s,\tau) = \sum_{k=0}^{N-1} \wdh\ell_k(s) \indi_{(t_k,t_{k+1}]}(\tau), \quad 0< \tau \leq s \leq T,
\eeq
and $\ell^\Dl(s,0):=\wdh\ell_0(s)$. Then, by the definitions of $P^\Dl$, $\eta^\Dl$, $\theta^\Dl$ and $\rho^\Dl$,
\beq{diag-Dl}
P^\Dl(s)=\Ga^\Dl(s,s),\quad \eta^\Dl(s)=\ga^\Dl(s,s),\quad \theta^\Dl(s)=\chi^\Dl(s,s),\quad \rho^\Dl(s)=\vartheta^\Dl(s,s),\quad s\in[0,T].
\eeq
\begin{lem}\label{lem:concat}
Let \assmref{ass} and \assmref{ass-rho} hold. For $0\leq k\leq N-1$, $(\wdh\Ga_k,\wdh\ga_k,\wdh\chi_k)$ is the mild solution on $[t_k,T]$ of \eqref{Ga-ga-chi-k}, with $(\Ga_k,\ga_k,\chi_k)$ replaced by $(\wdh\Ga_k,\wdh\ga_k,\wdh\chi_k)$ and with the same terminal conditions. Moreover, $\wdh\ga_k=\wdh\vartheta_kz^\Dl$ on $[t_k,T]$, $\wdh\vartheta_k$ satisfies \eqref{concat} on $[t_k,T]$, and
\beq{La-hat}
\wdh\Ga_kB+(\CC^\Dl)^*\wdh\Ga_kD-(\Theta^\Dl)^*R_k=\indi_{(t_{k+1},T]}\La_k\quad\text{on }(t_k,T].
\eeq
\end{lem}
\begin{proof}
On $(t_k,t_{k+1}]$ we have $\Theta^\Dl=\Theta_k=K_k^{-1}L_k$, $\Phi^\Dl=\Phi_k=-K_k^{-1}B^*\eta_k$ and $L_k^*=P_kB+C^*P_kD$. Hence
\bea
\ad P_k(A_1-B\Theta_k)+(A_1-B\Theta_k)^*P_k+\CC_k^*P_k\CC_k+\Theta_k^*R_k\Theta_k\\
\ad\quad=P_kA_1+A_1^*P_k+C^*P_kC-L_k^*\Theta_k-\Theta_k^*L_k+\Theta_k^*K_k\Theta_k\\
\ad\quad=P_kA_1+A_1^*P_k+C^*P_kC-L_k^*K_k^{-1}L_k,\\
\ad P_kB+\CC_k^*P_kD-\Theta_k^*R_k=L_k^*-\Theta_k^*K_k=0,\\
\ad \qv{P_kD\Phi_k,D\Phi_k}_H+2\qv{\eta_k,B\Phi_k}_H+\qv{R_k\Phi_k,\Phi_k}_U\\
\ad\quad=\qv{K_k\Phi_k,\Phi_k}_U+2\qv{B^*\eta_k,\Phi_k}_U=-\qv{K_k^{-1}B^*\eta_k,B^*\eta_k}_U.
\eea
So \eqref{P-k}, \eqref{eta-k} and \eqref{r-k} are the three equations of \eqref{Ga-ga-chi-k} on $[t_k,t_{k+1}]$ for $(P_k,\eta_k,\theta_k)$, and the second identity gives \eqref{La-hat}. Since $P_k(t_{k+1})=\Ga_k(t_{k+1})$, $\eta_k(t_{k+1})=\ga_k(t_{k+1})$ and $\theta_k(t_{k+1})=\chi_k(t_{k+1})$, the semigroup property gives the first claim. The relation $\wdh\ga_k=\wdh\vartheta_kz^\Dl$ follows from $\eta_k=\rho_kz_k$ and $\ga_k=\vartheta_kz^\Dl$, and \eqref{concat} was derived in the proof of \propref{prop:H3} from the mild forms of $\rho_k$ and $\vartheta_k$ only.
\end{proof}

We then obtain the following result.
\begin{prop}\label{prop:cauchy}
Let \assmref{ass} and \assmref{ass-rho} hold and $\{\Dl_\iota\}_{\iota=1}^\infty \subset \DD$, the set of all partitions of $[0,T]$, be such that $\Dl_{\iota+1}$ is a refinement of $\Dl_\iota,\forall\,\iota\in \NN$ and that $\lim_{ \iota \to \infty} \|\Dl_\iota\|=0$. Then, $\{\Ga^{\Dl_\iota}\}_{\iota=1}^\infty$, $\{\ga^{\Dl_\iota}\}_{\iota=1}^\infty$, and $\{\chi^{\Dl_\iota}\}_{\iota=1}^\infty$ are Cauchy sequences in $L^\infty(\Xi;\CL(H))$, $L^\infty(\Xi; H)$, and $L^\infty(\Xi; \rr)$, respectively.
\end{prop}
\begin{proof}
The proof follows the idea of \cite[Proposition 4.3]{DL20}, where the two-parameter functions are compared for fixed $\tau$. For $\iota,m\in \NN$, we consider $\Dl_\iota = \{0=t_0 < t_1< \cdots < t_N=T\}$ and $\Dl_{\iota+m} = \{0=\tilde t_{m_0} < \tilde t_{m_0+1} <\cdots < \tilde t_{m_1} < \cdots < \tilde t_{m_N}=T\}$ such that $\tilde t_{m_k}=t_k$ for $k=0,\cdots,N$. In what follows, $w\cd$ is a modulus of continuity, which may vary from line to line, and $\tau\in(\tilde t_i,\tilde t_{i+1}]\subset(t_k,t_{k+1}]$, so that $|\tilde t_i-t_k|\leq\|\Dl_\iota\|$. For $\tau=0$ we take $i=k=0$.

We first consider $\Ga^\Dl$. Set
\[
\nu_\Ga(s):=\sup_{0\leq\tau\leq s}|\Ga^{\Dl_{\iota+m}}(s,\tau)-\Ga^{\Dl_\iota}(s,\tau)|_{\CL(H)},\quad s\in[0,T].
\]
By \eqref{diag-Dl}, we have $|P^{\Dl_{\iota+m}}(r)-P^{\Dl_\iota}(r)|_{\CL(H)}\leq\nu_\Ga(r)$. Since $K^\Dl\geq\dl I$, the definitions of $\Theta^\Dl$ and $K^\Dl$ give
\beqa{diff-Theta-K}
\ad|\Theta^{\Dl_{\iota+m}}(r)-\Theta^{\Dl_\iota}(r)|+|(K^{\Dl_{\iota+m}})^{-1}(r)-(K^{\Dl_\iota})^{-1}(r)|\\
\aad\quad\leq c\,\big(\nu_\Ga(r)+|R^{\Dl_{\iota+m}}(r)-R^{\Dl_\iota}(r)|\big)\leq c\,\nu_\Ga(r)+w(\|\Dl_\iota\|).
\eeqa
By \lemref{lem:concat}, $\Ga^{\Dl_{\iota+m}}(\cdot,\tau)$ and $\Ga^{\Dl_\iota}(\cdot,\tau)$ are mild solutions on $[\tau,T]$ of the first equation of \eqref{Ga-ga-chi-k} with the weights at $\tilde t_i$ and $t_k$, respectively. By the uniform boundedness of $\Ga^\Dl$ and $\Theta^\Dl$, we have for $s\in[\tau,T]$,
\bea
\ad|\Ga^{\Dl_{\iota+m}}(s,\tau)-\Ga^{\Dl_\iota}(s,\tau)|\leq c\,\Big\{|(F+\bar F)(\tilde t_i)-(F+\bar F)(t_k)|+\int_s^T\big[\nu_\Ga(r)+|\Theta^{\Dl_{\iota+m}}(r)-\Theta^{\Dl_\iota}(r)|\\
\aad\quad+|(Q+\bar Q)(r,\tilde t_i)-(Q+\bar Q)(r,t_k)|+|R(r,\tilde t_i)-R(r,t_k)|\big]dr\Big\}.
\eea
By \eqref{diff-Theta-K} and the uniform continuity of $F,\bar F,Q,\bar Q$ and $R$, taking the supremum over $\tau$ gives $\nu_\Ga(s)\leq w(\|\Dl_\iota\|)+c\int_s^T\nu_\Ga(r)dr$. Gronwall's inequality yields
\beq{cauchy-Ga}
|\Ga^{\Dl_{\iota+m}}-\Ga^{\Dl_\iota}|_{L^\infty(\Xi;\CL(H))}\leq\sup_{s\in[0,T]}\nu_\Ga(s)\leq w(\|\Dl_\iota\|),
\eeq
and by \eqref{diff-Theta-K} the same bound holds for $|\Theta^{\Dl_{\iota+m}}(r)-\Theta^{\Dl_\iota}(r)|$ and $|(K^{\Dl_{\iota+m}})^{-1}(r)-(K^{\Dl_\iota})^{-1}(r)|$.

We next consider $\rho^\Dl$, which may jump at the nodes. Set
\[
\nu(s):=\sup_{0\leq\tau\leq s}|\vartheta^{\Dl_{\iota+m}}(s,\tau)-\vartheta^{\Dl_\iota}(s,\tau)|_{\CL(H)},\quad s\in[0,T],
\]
so that $|\rho^{\Dl_{\iota+m}}(s)-\rho^{\Dl_\iota}(s)|_{\CL(H)}\leq\nu(s)$ by \eqref{diag-Dl}. Recall that $\wdh\vartheta_k$ satisfies \eqref{concat} on $[t_k,T]$, where $\indi_{(t_{k+1},T]}\La_k$ is given by \eqref{La-hat}. By \eqref{concat} for both partitions and the uniform boundedness of $\Theta^\Dl$, $\rho^\Dl$, $\vartheta^\Dl$ and $\Ga^\Dl$, we have for $s\in[\tau,T]$,
\beqa{diff-rho}
\ad|\vartheta^{\Dl_{\iota+m}}(s,\tau)-\vartheta^{\Dl_\iota}(s,\tau)|\leq c\,\Big\{|[\bar F\bar G](\tilde t_i)-[\bar F\bar G](t_k)|+\int_s^T\big[\nu(r)+\nu_\Ga(r)\\
\aad\quad+|\Theta^{\Dl_{\iota+m}}(r)-\Theta^{\Dl_\iota}(r)|+|(K^{\Dl_{\iota+m}})^{-1}(r)-(K^{\Dl_\iota})^{-1}(r)|\\
\aad\quad+|R(r,\tilde t_i)-R(r,t_k)|+|[\bar Q\bar S](r,\tilde t_i)-[\bar Q\bar S](r,t_k)|\big]dr\Big\}.
\eeqa
By \eqref{cauchy-Ga} and the uniform continuity of $\bar F\bar G$, $\bar Q\bar S$ and $R$, taking the supremum over $\tau$ gives $\nu(s)\leq w(\|\Dl_\iota\|)+c\int_s^T\nu(r)dr$, and Gronwall's inequality yields
\beq{cauchy-rho}
|\rho^{\Dl_{\iota+m}}-\rho^{\Dl_\iota}|_{L^\infty([0,T];\CL(H))}\leq\sup_{s\in[0,T]}\nu(s)\leq w(\|\Dl_\iota\|)\to0\ \text{ as }\iota\to\infty.
\eeq
Consequently, $|\Psi^{\Dl_{\iota+m}}(r)-\Psi^{\Dl_\iota}(r)|\leq w(\|\Dl_\iota\|)$.

From \eqref{rho-z-Dl}, the equation governing $z^\Dl$ is linear with respect to $\rho^\Dl$, and $z^\Dl$ is continuous on $[0,T]$. An analogous argument establishes that
\beq{cauchy-z}
|z^{\Dl_{\iota+m}}-z^{\Dl_\iota}|_{L^\infty([0,T];H)}=\sup_{s\in [0,T]}
|z^{\Dl_{\iota+m}}(s)-z^{\Dl_\iota}(s)|_H \leq w( \|\Dl_\iota\|) \to 0.
\eeq

Finally, by \lemref{lem:concat}, $\ga^\Dl(s,\tau)=\vartheta^\Dl(s,\tau)z^\Dl(s)$. Hence, by \eqref{cauchy-rho} and \eqref{cauchy-z},
\[
|\ga^{\Dl_{\iota+m}}- \ga^{\Dl_{\iota}}|_{L^\infty(\Xi; H)} \leq c\sup_{s\in[0,T]}\big(\nu(s)+|z^{\Dl_{\iota+m}}(s)-z^{\Dl_\iota}(s)|_H\big)\leq w(\|\Dl_\iota\|).
\]
Since $\Phi^\Dl=\Psi^\Dl z^\Dl$, the Cauchy property of $\{\chi^{\Dl_\iota}\}_{\iota=1}^\infty$ follows in the same way from the third equation of \eqref{Ga-ga-chi-k} for $\wdh\chi_k$ and that of $\{\Ga^{\Dl_\iota}\}_{\iota=1}^\infty$, $\{\ga^{\Dl_\iota}\}_{\iota=1}^\infty$, $\{z^{\Dl_\iota}\}_{\iota=1}^\infty$ and $\{\Psi^{\Dl_\iota}\}_{\iota=1}^\infty$. This completes the proof.
\end{proof}

By the definitions of $\Theta^\Delta$ and $\Phi^\Delta$ in \eqref{Theta-Phi-Dl}, analyzing their convergence requires a prior investigation into the convergence of $P^\Delta$ and $\eta^\Delta$, which in turn depends on that of $\Gamma^\Delta$ and $\gamma^\Delta$. To this end, we establish the following convergence theorem.
\begin{theorem}\label{thm:conv}
Suppose that \assmref{ass} and \assmref{ass-rho} hold. Then there are unique mild solutions $\Ga(\cdot,\cdot), \ga(\cdot,\cdot)$ and $\chi(\cdot,\cdot)$. Moreover, the following convergence holds
\beq{conv-tau-s}
\lim_{\|\Dl\| \to 0} \big(|\Ga^\Dl(s,\tau)-\Ga(s,\tau)|_{\CL(H)}+|\ga^\Dl(s,\tau)-\ga(s,\tau)|_H + |\chi^\Dl(s,\tau)-\chi(s,\tau)|\big)=0
\eeq
uniformly in $(s,\tau)\in \Xi$. Furthermore, it holds that
\beq{conv-ss}
\lim_{\|\Dl\|\to 0} \big( |P^\Dl(s)-\Ga(s,s)|_{\CL(H)}+ |\eta^\Dl(s)-\ga(s,s)|_H + |\theta^\Dl(s)-\chi(s,s)| \big) = 0,
\eeq
which implies that
\bea\ad
\lim_{\|\Dl\|\to 0} (|\Theta^\Dl(s)-\Theta(s)|_{\CL(H;U)}+|\Phi^\Dl(s)-\Phi(s)|_U )=0.
\eea
where $\Theta(s)$ and $\Phi(s)$ are defined in \eqref{def-theta} and \eqref{def-Phi}, respectively.
\end{theorem}
\begin{proof}
\textbf{Step 1 (Existence of mild solutions)}. Let $\{\Dl_\iota\}_{\iota=1}^\infty$ be a sequence of partitions of $[0,T]$ such that $\Dl_{\iota+1}$ is a refinement of $\Dl_\iota$ for all $\iota\in \NN$ and that $\lim_{\iota\to \infty}\|\Dl_\iota\|=0$. By \propref{prop:cauchy}, $\{\Ga^{\Dl_\iota}, \ga^{\Dl_\iota}, \chi^{\Dl_\iota}\}_{\iota=1}^\infty$ is a Cauchy sequence in $L^\infty(\Xi; \CL(H) \times H\times \rr)$. Therefore, there exist $(\Ga,\ga,\chi) \in L^\infty(\Xi; \CL(H) \times H \times \rr)$ such that
\beqa{conv-Ga-Dl}\ad\lim_{\iota\to \infty} \big( |\Ga^{\Dl_\iota}(\cdot,\cdot)- \Ga(\cdot,\cdot)|_{L^\infty(\Xi; \CL(H))}+|\ga^{\Dl_\iota}(\cdot,\cdot)- \ga(\cdot,\cdot)|_{L^\infty(\Xi; H)} \\
\aad\qquad\quad + |\chi^{\Dl_\iota}(\cdot,\cdot)- \chi(\cdot,\cdot)|\big) =0.
\eeqa
By the estimates in the proof of \propref{prop:cauchy}, the sequences $\{\Ga^{\Dl_\iota}\}$, $\{\ga^{\Dl_\iota}\}$, $\{\chi^{\Dl_\iota}\}$ and $\{\vartheta^{\Dl_\iota}\}$ are uniformly Cauchy at every point of $\Xi$. So the limits are defined pointwise on $\Xi$, in particular on the diagonal $\tau=s$. Moreover, these estimates hold for any partition $\Dl$ and any refinement $\Dl'$ of $\Dl$ with the same modulus $w(\|\Dl\|)$. Comparing two partitions through their common refinement, we see that the limits do not depend on the sequence $\{\Dl_\iota\}$ and that \eqref{conv-tau-s} holds as $\|\Dl\|\to0$.

By \eqref{diag-Dl}, we have $P^{\Dl_\iota}(s)=\Ga^{\Dl_\iota}(s,s)$, $\eta^{\Dl_\iota}(s)=\ga^{\Dl_\iota}(s,s)$ and $\theta^{\Dl_\iota}(s)=\chi^{\Dl_\iota}(s,s)$. Hence the uniform convergence at every point of $\Xi$ gives
\beqa{conv-P-Dl}\ad
\lim_{\iota \to \infty}\sup_{s\in [0,T]}\big(|P^{\Dl_\iota}(s)-\Ga(s,s)|_{\CL(H)}+ |\eta^{\Dl_\iota}(s)-\ga(s,s)|_H + |\theta^{\Dl_\iota}(s)-\chi (s,s)|\big) =0.
\eeqa
Consequently, it follows that
\bea\ad
\lim_{\iota\to \infty}(|\Theta^{\Dl_\iota}(s)-\Theta(s)|_{\CL(H;U)}+|\Phi^{\Dl_\iota}(s)-\Phi(s)|_U )=0,
\eea
where
\bea
\Theta(s)\ad=[R(s,s)+D^{*} \Ga(s,s)D]^{-1}[B^*\Ga(s,s)+D^* \Ga(s,s)C] \\
\Phi(s)\ad= -[R(s,s)+D^{*} \Ga(s,s) D]^{-1} B^* \ga(s,s).
\eea
Moreover, by the proof of \propref{prop:cauchy}, $\rho^{\Dl_\iota}$ and $z^{\Dl_\iota}$ also converge. We have
\[
\lim_{\iota\to \infty}\sup_{s\in [0,T]}|z^{\Dl_\iota}(s)-z(s)|=0,
\]
where $z$ satisfies the second equation in \eqref{fb-eq}.

We now show that $\Ga,\ga$, and $\chi$ are the mild solutions of \eqref{bRic} and \eqref{fb-eq}. For any $\phi\in H$, from \lemref{lem:concat}, we have for $s\in[t_k,T]$
\bea
\wdh\Ga_{k}(s)\phi \ad\!\!\!= \Big\{e^{(T-s)A^*} (F(t_k)+\bar F(t_k))e^{(T-s)A} + \int_s^T e^{(r-s)A^*} \big[\wdh\Ga_{k}(r) (A_1 - B\Theta^{\Dl_\iota})\\
\aad\quad + (A_1-B\Theta^{\Dl_\iota})^* \wdh\Ga_{k}(r) + Q(r,t_k)+\bar Q(r,t_k) \\
\aad\quad + (C-D\Theta^{\Dl_\iota})^* \wdh\Ga_{k}(r) (C-D\Theta^{\Dl_\iota})+ (\Theta^{\Dl_\iota})^* R(r,t_k) \Theta^{\Dl_\iota} \big] e^{(r-s)A}  dr \Big\} \phi.
\eea
Together with the definition of $\Ga^{\Dl_\iota}$ in \eqref{def-ell-Dl-s}, this gives for any $\tau\in (t_k,t_{k+1}]$ and $s\in[\tau,T]$,
\beqa{Ga-Dl-taus}\ad\!\!\!
\Ga^{\Dl_\iota}(s,\tau)\phi  =\Big\{e^{(T-s)A^*}(F(t_k)+\bar F(t_k))e^{(T-s)A} +\int_s^T e^{(r-s)A^*}\big[ \Ga^{\Dl_\iota}(r,\tau)\\
\aad\quad \times (A_1- B \Theta^{\Dl_\iota})+(A_1 - B \Theta^{\Dl_\iota})^* \Ga^{\Dl_\iota}(r,\tau) + Q(r,t_k)+ \bar Q(r,t_k)\\
\aad\quad + (C-D\Theta^{\Dl_\iota})^* \Ga^{\Dl_\iota}(r,\tau) (C-D\Theta^{\Dl_\iota})+ (\Theta^{\Dl_\iota})^* R(r,t_k) \Theta^{\Dl_\iota}
\big]e^{(r-s)A} dr \Big\} \phi.
\eeqa

The bound $t_{k+1} - t_k\leq \|\Dl_\iota\|\to0$ implies that $t_k \to \tau$ as $\iota \to \infty$. From \eqref{conv-Ga-Dl}--\eqref{Ga-Dl-taus}, it follows that as $\iota \to \infty$,
\beqa{Ga-tau}
\Ga(s,\tau)\phi \ad = \Big\{ e^{(T-s)A^*}(F(\tau)+\bar F(\tau)) e^{(T-s)A} + \int_s^T e^{(r-s)A^*} \big[\Ga(r,\tau)(A_1-B \Theta)\\
\aad\quad + (A_{1}-B\Theta)^* \Ga(r,\tau) + Q(r,\tau)+\bar Q(r,\tau) + (C-D\Theta)^* \Ga(r,\tau) (C-D\Theta) \\
\aad\quad + (\Theta)^* R(r,\tau) \Theta
\big] e^{(r-s)A} dr \Big\}\phi.
\eeqa
Consequently, $\Ga(\cdot,\cdot)$ is strongly continuous with respect to the first variable.

From the second equation of \eqref{Ga-ga-chi-k} and \lemref{lem:concat}, for any $\tau\in (t_k,t_{k+1}]$ and $s\in[\tau,T]$, we have
\bea\!\!\!
\wdh\ga_{k}(s)\ad\!\!\!\!= e^{(T-s)A^*}(-[\bar F \bar G](t_k) z^{\Dl_\iota}(T))\!+\!\! \int_s^T\! e^{(r-s)A^*} \big[ (A_1 -B \Theta^{\Dl_\iota})^{*}\wdh\ga_{k}(r)+ \wdh\Ga_{k}(r) B \Phi^{\Dl_\iota} \\
\aad\!\!\! + (C-D\Theta^{\Dl_\iota})^* \wdh\Ga_{k}(r) D \Phi^{\Dl_\iota} - (\Theta^{\Dl_\iota})^* R(r,t_k)  \Phi^{\Dl_\iota} - \bar Q(r,t_k) \bar S(r,t_k) z^{\Dl_\iota}(r) \big] dr \\
\ad\!\!\! = e^{(T-s)A^*}(-[\bar F \bar G](t_k) z^{\Dl_\iota}(T)) +\int_s^T e^{(r-s)A^*} \big[(A_1-B\Theta^{\Dl_\iota}(r))^{*} \ga^{\Dl_\iota}(r,\tau) \\
\aad\!\!\! + \Ga^{\Dl_\iota}(r,\tau) B \Phi^{\Dl_\iota} + (C-D \Theta^{\Dl_\iota})^* \Ga^{\Dl_\iota}(r,\tau) D \Phi^{\Dl_\iota} \\
\aad\!\!\! - (\Theta^{\Dl_\iota})^* R(r,t_k)  \Phi^{\Dl_\iota} - \bar Q(r,t_k) \bar S(r,t_k) z^{\Dl_\iota}(r) \big] dr,
\eea
where $\wdh\ga_k(s)=\ga^{\Dl_\iota}(s,\tau)$ and
\bea
z^{\Dl_\iota}(s)= e^{sA}z_0 + \int_0^s e^{(s-r)A} [(A_1-B\Theta^{\Dl_\iota}(r)) z^{\Dl_\iota}(r)+ B \Phi^{\Dl_\iota}(r)] dr.
\eea
Letting $\iota \to \infty$, from  \eqref{conv-Ga-Dl} and \eqref{conv-P-Dl}, we have
\beqa{ga-tau}
\ga(s,\tau)&\!\!\!= e^{(T-s)A^*}(- [\bar F \bar G](\tau) z(T)) + \int_s^T e^{(r-s)A^*} \big[(A_1-B\Theta)^{*}\ga(r,\tau)+\Ga(r,\tau) B \Phi \\
& + (C- D\Theta)^* \Ga(r,\tau) D \Phi - \Theta^* R(r,\tau) \Phi - \bar Q(r,\tau) \bar S(r,\tau) z(r) \big]dr,
\eeqa
which is coupled with
\beqa{z-tau}
z(s) = e^{sA}z_0 + \int_0^s e^{(s-r)A} \big[(A_1-B \Theta) z(r) + B \Phi(r)\big]dr.
\eeqa

Finally, repeating this procedure, it is not difficult to obtain that
\beqa{chi-tau}
\chi(s,\tau)& \!\!\!=\qv{\bar F(\tau) \bar G(\tau)z(T), \bar G(\tau)z(T)}_H + \int_s^T \qv{\Ga(r,\tau) D\Phi, D\Phi}_H  \\
& \!\!\!\! + 2\qv{\ga(r,\tau), B \Phi}_H +\qv{R(r,\tau)\Phi, \Phi}_U + \qv{\bar Q(r,\tau)\bar S(r,\tau) z(r), \bar S(r,\tau)z(r)} dr.
\eeqa

\smallskip
\noindent
\textbf{Step 2 (Uniqueness of the solutions)}. The uniqueness of $\Gamma$ can be found in \cite[pp. 505--507]{DL20}; the details are omitted here for brevity. In particular, Dou and L\"{u} \cite[p. 507]{DL20} proved that
\beq{est-til-Ga}
|\tilde \Ga(s,\tau)|_{\CL(H)}=0,\quad \text{ for } (s,\tau) \in [0,T] \times [0,s],
\eeq
where $\tilde\Ga:=\Ga^1-\Ga^2$ with  $\Ga^1,\Ga^2$ being two solutions of \eqref{Ga-tau}.

We proceed to prove the uniqueness of $(\ga,z)$. Let $(\ga^1,z^1),(\ga^2,z^2)$ be two mild solutions to \eqref{ga-tau} and \eqref{z-tau}. Define $\tilde \ga:=\ga^1-\ga^2$ and $\tilde z:= z^1-z^2$. Then we have
\beqa{til-ga}
\tilde\ga(s,\tau)\ad\!\!\! = e^{(T-s)A^*}(-[\bar F \bar G](\tau) \tilde z(T)) + \int_s^T e^{(r-s)A^*} \big[(A_1-B\Theta^1)^{*}\tilde\ga(r,\tau)+ \tilde\Ga(r,\tau) B\Phi^1 \\
\aad\quad + (C-D\Theta^2)^* \tilde \Ga(r,\tau) D\Phi^1  -\bar Q(r,\tau)\bar S(r,\tau) \tilde z(r) \big]dr \\
\aad - \int_s^T e^{(r-s)A^*}\big[(\Theta^1-\Theta^2)^*B^* \ga^2(r,\tau) - \Ga^2(r,\tau)B (\Phi^1-\Phi^2) \\
\aad\quad +(D\Theta^1-D\Theta^2)^* \Ga^1(r,\tau) D \Phi^1 - (C-D\Theta^2)^* \Ga^2(r,\tau)(D\Phi^1-D\Phi^2)
\\
\aad\quad +(\Theta^1)^*R(r,\tau) \Phi^1-(\Theta^2)^* R(r,\tau) \Phi^2 \big] dr,
\eeqa
and
\beqa{til-z}
\tilde z(s)=\int_0^s e^{(s-r)A} \big[(A_1-B\Theta^1)\tilde z(r)- B(\Theta^1-\Theta^2)z^2(r)+ B\Phi^1-B\Phi^2 \big]dr,
\eeqa
where $\Theta^j$ and $\Phi^j, j=1,2$ are given by
\beqa{Theta-Phi-j}
\Theta^j(r)\ad = [R(r,r)+D^*\Ga^j(r,r) D]^{-1} [B^* \Ga^j(r,r)+D^*\Ga^j(r,r)C], \\
\Phi^j(r) \ad = -[R(r,r)+D^* \Ga^j(r,r) D]^{-1} B^* \ga^j(r,r).
\eeqa

To obtain the desired estimates for \eqref{til-ga} and \eqref{til-z}, we first need to estimate $\Theta^1 - \Theta^2$ and $\Phi^1 - \Phi^2$. From \eqref{Theta-Phi-j}, it follows that
\bea\ad
|\Theta^1(r)-\Theta^2(r)|_{\CL(H;U)}\\
\aad \leq \big|[R(r,r)+D^*\Ga^1(r,r)D]^{-1}-[R(r,r)+D^* \Ga^2(r,r)D]^{-1} \big|_{\CL(U)}\\
\aad\quad \times \big|B^* \Ga^2(r,r)+D^* \Ga^2(r,r)C \big|_{\CL(H;U)} \\
\aad+ \big|[R(r,r)+D^* \Ga^1(r,r)D]^{-1}\\
\aad\quad \times \big[ B^* (\Ga^1(r,r)-\Ga^2(r,r))+ D^*(\Ga^1(r,r)-\Ga^2(r,r))C \big] \big|_{\CL(H;U)}
\eea
In the above, by the local Lipschitz
property, there exists a constant $c(\Ga^1,\Ga^2)>0$ s.t.
\bea\ad
\big|[R(r,r)+D^* \Ga^1(r,r) D]^{-1} -[R(r,r)+D^* \Ga^2(r,r)D]^{-1} \big|_{\CL(U)} \\
\aad \leq \bigg( \prod_{j=1}^2 |[R(r,r)+D^* \Ga^j(r,r)D]^{-1}|_{\CL(U)}\bigg)  \times |D^*(\Ga^1(r,r)-\Ga^2(r,r))D|_{\CL(U)} \\
\aad \leq c(\Ga^1,\Ga^2)|\Ga^1(r,r)-\Ga^2(r,r)|_{\CL(H)},
\eea
which yields that
\beqa{Theta-12}
|\Theta^1(r)-\Theta^2(r)|_{\CL(H;U)} \leq c(\Ga^1,\Ga^2)\,|\Ga^1(r,r)-\Ga^2(r,r)|_{\CL(H)}.
\eeqa
Analogously, there exists a constant
$c(\Ga^1,\Ga^2,\ga^1,\ga^2)>0$ such that
\beqa{Phi-12}\ad
|\Phi^{1}(r)-\Phi^{2}(r)|_{U} \\
\aad \leq \big|[R(r,r)+D^*\Ga^1(r,r)D]^{-1}-[R(r,r)+D^* \Ga^2(r,r)D]^{-1} \big|_{\CL(U)} |B^{*}\ga^2(r,r)|_{U} \\
\aad\quad  +|[R(r,r)+D^* \Ga^1(r,r)D]^{-1} B^{*}[\ga^1(r,r)-\ga^2(r,r)]|_{U} \\
\aad\leq c(\Ga^1,\Ga^2,\ga^1,\ga^2)\big\{|\Ga^1(r,r)-\Ga^2(r,r)|_{\CL(H)}+|\ga^1(r,r)-\ga^2(r,r)|_H \big\}.
\eeqa
Here, $c(\Ga^1,\Ga^2), c(\Ga^1,\Ga^2,\ga^1,\ga^2)$ are constants depending on $\Ga^1,\Ga^2,\ga^1$ and $\ga^2$, which may change from line to line.

For any fixed $\tau\in [0,T]$ and $s\in[\tau,T]$, equations \eqref{til-ga} and \eqref{til-z} can be written in the following differential form
\beqa{tilde-ga-z}\left\{\barray
- \partial_s \tilde \ga(s,\tau) =[A+A_1-B\Theta^{1}]^*\tilde \ga(s,\tau)+ \tilde \Ga(s,\tau) B \Phi^1 + (C-D\Theta^2)^* \tilde \Ga(s,\tau) D\Phi^1\\
\qquad\qquad\quad - \bar Q(s,\tau) \bar S(s,\tau)\tilde z(s) -(\Theta^1-\Theta^2)^*B^*\ga^2(s,\tau)+\Ga^2(s,\tau) B(\Phi^1-\Phi^2) \\
\qquad\qquad\quad - (D\Theta^1-D\Theta^2)^* \Ga^1(s,\tau) D\Phi^1+(C-D\Theta^2)^* \Ga^2(s,\tau) (D\Phi^1-D\Phi^2)\\
\qquad\qquad\quad -(\Theta^1)^* R(s,\tau)  \Phi^1+ (\Theta^2)^* R(s,\tau)\Phi^2  \\
\partial_s \tilde z(s) =[A+A_1-B\Theta^1] \tilde z(s)- B(\Theta^1-\Theta^2) z^2(s)+B(\Phi^1-\Phi^2) \\
\tilde \ga(T,\tau) = -[\bar F \bar G](\tau) \tilde z(T), \quad \tilde z(0)=0.
\earray\right.\eeqa

Define $\tilde\Phi:=\Phi^1-\Phi^2$. By \eqref{est-til-Ga}, $\Ga^1=\Ga^2$. Hence $\Theta^1=\Theta^2=:\Theta$, and $K(s):=R(s,s)+D^*\Ga^1(s,s)D$ is the same for both solutions. Thus \eqref{tilde-ga-z} reduces to
\beq{tilde-ga-z-red}\left\{\barray
-\partial_s\tilde\ga(s,\tau)=[A+A_1-B\Theta]^*\tilde\ga(s,\tau)+\La(s,\tau)\tilde\Phi(s)-\bar Q(s,\tau)\bar S(s,\tau)\tilde z(s),\\
\partial_s\tilde z(s)=[A+A_1-B\Theta]\tilde z(s)+B\tilde\Phi(s),\\
\tilde\ga(T,\tau)=-[\bar F\bar G](\tau)\tilde z(T),\quad \tilde z(0)=0,
\earray\right.\eeq
where $\La(s,\tau):=\Ga^1(s,\tau)B+(C-D\Theta)^*\Ga^1(s,\tau)D-\Theta^*R(s,\tau)$ and $\tilde\Phi(s)=-K(s)^{-1}B^*\tilde\ga(s,s)$.

Let $\vartheta^\Dl(s,\tau)$ be defined by \eqref{def-ell-Dl-s}. By the proof of \propref{prop:cauchy}, $\{\vartheta^{\Dl_\iota}\}_{\iota=1}^\infty$ is a Cauchy sequence in $L^\infty(\Xi;\CL(H))$ and $\{\Psi^{\Dl_\iota}\}_{\iota=1}^\infty$ converges uniformly. Denote the limits by $\vartheta$ and $\Psi$. Since $\rho^\Dl(s)=\vartheta^\Dl(s,s)$ by \eqref{diag-Dl}, we have $\Psi(s)=-K(s)^{-1}B^*\vartheta(s,s)$. Letting $\iota\to\infty$ in \eqref{concat} and using \eqref{La-hat}, $\vartheta$ is the mild solution of
\beq{vartheta-lim}\left\{\barray
-\partial_s\vartheta(s,\tau)=[A+A_1-B\Theta]^*\vartheta(s,\tau)+\vartheta(s,\tau)[A+A_1-B\Theta+B\Psi(s)]\\
\qquad\qquad\qquad\;\, +\La(s,\tau)\Psi(s)-\bar Q(s,\tau)\bar S(s,\tau),\\
\vartheta(T,\tau)=-[\bar F\bar G](\tau).
\earray\right.\eeq
Now let $\mu(s,\tau):=\tilde\ga(s,\tau)-\vartheta(s,\tau)\tilde z(s)$. Since $\tilde\Phi(s)=\Psi(s)\tilde z(s)-K(s)^{-1}B^*\mu(s,s)$, it follows from \eqref{tilde-ga-z-red} and \eqref{vartheta-lim} that
\bea
-\partial_s\mu(s,\tau)\ad=-\partial_s\tilde\ga(s,\tau)+\partial_s\vartheta(s,\tau)\tilde z(s)+\vartheta(s,\tau)\partial_s\tilde z(s)\\
\ad=[A+A_1-B\Theta]^*\mu(s,\tau)+[\La(s,\tau)+\vartheta(s,\tau)B]\big(\tilde\Phi(s)-\Psi(s)\tilde z(s)\big)\\
\ad=[A+A_1-B\Theta]^*\mu(s,\tau)-[\La(s,\tau)+\vartheta(s,\tau)B]K(s)^{-1}B^*\mu(s,s),
\eea
with $\mu(T,\tau)=0$. As in \thmref{thm:unique}, this computation is justified by the Yosida approximation, and $\mu$ has the mild form
\bea
\mu(s,\tau)=\int_s^T e^{(r-s)A^*}\big\{(A_1-B\Theta)^*\mu(r,\tau)-[\La(r,\tau)+\vartheta(r,\tau)B]K(r)^{-1}B^*\mu(r,r)\big\}dr.
\eea
Let $\nu(s):=\sup_{0\leq \tau\leq s}|\mu(s,\tau)|_H$. Since $|\mu(r,\tau)|_H\vee|\mu(r,r)|_H\leq \nu(r)$ for $\tau\leq s\leq r$, we obtain $\nu(s)\leq c\int_s^T \nu(r)dr$, and Gronwall's inequality gives $\nu\equiv0$. Hence $\tilde\ga(s,\tau)=\vartheta(s,\tau)\tilde z(s)$ and $\tilde\Phi=\Psi\tilde z$. Then $\partial_s\tilde z(s)=[A+A_1-B\Theta+B\Psi(s)]\tilde z(s)$ with $\tilde z(0)=0$, which yields $\tilde z\equiv 0$ and then $\tilde\ga\equiv0$.

Hence $\ga^1=\ga^2$ and $z^1=z^2$, which establishes the uniqueness of \eqref{ga-tau} and \eqref{z-tau}. The uniqueness of $\chi$ in \eqref{chi-tau} follows that of $(\Ga, \ga, z)$. The proof is completed.
\end{proof}

\section{Existence of Closed-loop Equilibrium Strategy}\label{sec:exist}
In this section, we establish that $\Theta,\Phi$ obtained in the previous section are a closed-loop equilibrium strategy.

\begin{prop}\label{prop:equi}
Suppose that \assmref{ass} and \assmref{ass-rho} hold. For $x\in L^2_{\CF_t}(\Omega;H)$, let $\wdh X$ and $\wdh u$ be the equilibrium state and equilibrium control given by \eqref{sys-fed} and \eqref{u-hat}. Then for any $u\cd \in \CU[t,T]$,  we have
\beq{J-inf}
\liminf_{\e\to 0^+} \frac{J(t,\wdh X(t); u^\e\cd)-J(t,\wdh X(t);\wdh u\cd)}{\e} \geq 0, \quad \PP\text{-a.s.},
\eeq
where $u^\e\cd=u\cd \indi_{[t,t+\e)}\cd + [-\Theta\cd X^\e\cd + \Phi\cd] \indi_{[t+\e,T]}\cd$, with $X^\e\cd$ being the solution of \eqref{see} corresponding to the control $u^\e \cd$.
\end{prop}
\begin{proof}

Fix $t\in[0,T)$ and set $x:=\wdh X(t)$. For $s\in[t,T]$ and $y\in H$, let
\[
V^t(s,y):=\half\big[\qv{\Ga(s,t)y,y}_H+2\qv{\ga(s,t),y}_H+\chi(s,t)\big].
\]
By \eqref{Ga-tau}, \eqref{ga-tau} and \eqref{chi-tau}, applying It\^o's formula to $V^t(s,X(s))$ along the closed-loop system with the feedback $(\Theta,\Phi)$, which is justified by the Yosida approximation as in the proof of \thmref{thm:unique}, shows that $V^t(s,y)$ is the cost on $[s,T]$ with the weights at time $t$ when $(\Theta,\Phi)$ is used on $[s,T]$ from the state $y$. Hence $J(t,x;\wdh u)=V^t(t,x)$ and
\[
J(t,x;u^\e)=\EE_t\Big[\int_t^{t+\e}f(t,s,X^\e(s),\pi(s),u(s))ds+V^t(t+\e,X^\e(t+\e))\Big].
\]
Applying It\^o's formula to $V^t(s,X^\e(s))$ on $[t,t+\e]$ and using the equations of $\Ga(\cdot,t)$, $\ga(\cdot,t)$ and $\chi(\cdot,t)$, all terms that do not depend on the control cancel, and we obtain
\beq{J-diff}
J(t,x;u^\e)-J(t,x;\wdh u)=\EE_t\int_t^{t+\e}\big[\I(s,X^\e(s),u(s))-\I(s,X^\e(s),\hat v(s))\big]ds,
\eeq
where $\hat v(s):=-\Theta(s)X^\e(s)+\Phi(s)$ and
\[
\I(s,y,v):=\half\qv{R(s,t)v,v}_U+\qv{Bv,\Ga(s,t)y+\ga(s,t)}_H+\half\qv{\Ga(s,t)(Cy+Dv),Cy+Dv}_H.
\]
Let $K(s,t):=R(s,t)+D^*\Ga(s,t)D$ and $L(s,t):=B^*\Ga(s,t)+D^*\Ga(s,t)C$. By \eqref{def-theta} and \eqref{def-Phi}, $K(s,s)\hat v(s)+L(s,s)X^\e(s)+B^*\ga(s,s)=0$. Hence, with
\[
\kappa(s):=[K(s,t)-K(s,s)]\hat v(s)+[L(s,t)-L(s,s)]X^\e(s)+B^*[\ga(s,t)-\ga(s,s)],
\]
and since $K(s,t)\geq\dl I$, we have
\bea
\I(s,X^\e,u)-\I(s,X^\e,\hat v)\ad=\half\qv{K(s,t)(u-\hat v),u-\hat v}_U+\qv{\kappa(s),u-\hat v}_U\\
\ad\geq\frac{\dl}{2}|u-\hat v|_U^2-|\kappa(s)|_U|u-\hat v|_U\geq-\frac{1}{2\dl}|\kappa(s)|_U^2.
\eea
By \eqref{Ga-tau}, \eqref{ga-tau}, Gronwall's inequality and the uniform continuity of $F,\bar F,\bar G,Q,\bar Q,\bar S$ and $R$, there is a modulus of continuity $w$ such that $|\Ga(s,t)-\Ga(s,s)|_{\CL(H)}+|\ga(s,t)-\ga(s,s)|_H\leq w(s-t)$. Together with the boundedness of $\Theta$ and $\Phi$, this gives $|\kappa(s)|_U\leq w(\e)(1+|X^\e(s)|_H)$ for $s\in[t,t+\e]$. Thus, by \eqref{J-diff},
\bea
J(t,x;u^\e)-J(t,x;\wdh u)\ad\geq-\frac{w(\e)^2}{\dl}\,\EE_t\int_t^{t+\e}\big(1+|X^\e(s)|_H^2\big)ds\\
\ad\geq-c\,\e\, w(\e)^2\Big(1+|x|_H^2+\EE_t\int_t^T|u(s)|_U^2ds\Big),
\eea
which implies \eqref{J-inf}.
\end{proof}
Combining \thmref{thm:conv} and \propref{prop:equi}, we obtain \thmref{thm:main}.

\end{document}